\documentclass[a4paper,11pt]{article}
\usepackage{amsmath,amsfonts,amssymb,amsthm}
\usepackage{graphicx,subfigure,booktabs,multirow}
\usepackage{mathrsfs,enumerate,color,bm}%,mathdots
\usepackage{algorithm,algpseudocode}
\usepackage[textwidth=14cm]{geometry}
\usepackage{color}
\usepackage{hyperref}

\newtheorem{theorem}{Theorem}
\newtheorem{lemma}[theorem]{Lemma}
\newtheorem{proposition}{Proposition}

\newtheorem{definition}{Definition}

\newcommand{\rank}{\operatorname{rank}}

\newcommand{\diag}{\operatorname{diag}}

\graphicspath{}

\title{Nonnegative Low Rank Matrix Approximation for Nonnegative Matrices}

\author{Guang-Jing Song\thanks{School of Mathematics and Information Sciences, Weifang University,
Weifang 261061, P.R. China. (email: sgjshu@163.com)}
\and
Michael K. Ng\thanks{Department of Mathematics, The University of Hong Kong, Pokfulam, Hong Kong
(email: mng@maths.hku.hk).
M. Ng's research supported in part by the HKRGC GRF 12306616,
12200317, 12300218 and 12300519, and HKU 104005583.}
}

\begin{document}
\maketitle

\begin{abstract}
This paper describes a new algorithm for computing Nonnegative Low Rank Matrix (NLRM) approximation for
nonnegative matrices. Our approach is completely different from classical nonnegative matrix
factorization (NMF) which has been studied for more than twenty five years. For a given nonnegative matrix, the usual NMF
approach is to determine two nonnegative low rank matrices such that the distance between their product and the
given nonnegative matrix is
as small as possible. However, the proposed NLRM approach is to determine a nonnegative low rank matrix such that
the distance between such matrix and the given nonnegative matrix is as small as possible.
There are two advantages. (i) The minimized distance by the proposed NLRM method can be smaller than
that by the NMF method, and it implies that the proposed NLRM method can obtain a better low rank matrix
approximation. (ii) Our low rank matrix admits a matrix singular value decomposition automatically
which provides a significant index based on singular values that can be used to
identify important singular basis vectors, while this information cannot be
obtained in the classical NMF. The proposed NLRM approximation
algorithm was derived using the alternating projection on the low rank matrix manifold and the non-negativity
property.
Experimental results are presented to demonstrate the above mentioned advantages of the proposed NLRM method
compared the NMF method.
\end{abstract}

\vspace{3mm}
\noindent
Keywords: Nonnegative matrix, low-rank approximation, manifolds, projections

\vspace{3mm}
\noindent
AMS subject classiflcations. 15A23, 65f22.

\section{Introduction}

Nonnegative data matrices appear in many data analysis applications.
For instance, in image analysis, image pixel values are nonnegative and the associated nonnegative image data matrices can be
formed for clustering and recognition \cite{Chen13,Ding05,Ding06,Guil02,Guil03,Jing12,Lee99,Liu12,Liu15,Wang05,Zhan05}.
In text mining, the frequencies of terms in documents are nonnegative and the resulted nonnegative
term-to-document data matrices can be constructed for clustering \cite{Berr10,Li06,Pauc04,Xu03}.
In bioinformatics, nonnegative gene expression values are studied and nonnegative gene expression data matrices
are generated for diseases and genes classification \cite{ Cho04,Cich09,Gao06,Kim03,Kim07,Pasc06,Wang06}.
Nonnegative Low Rank Matrix (NLRM) approximation for nonnegative matrices
play a key role in all these applications. Its main purpose is to identify a latent feature
space for objects representation. The classification, clustering or recognition analysis can be
done by using these latent features. Lee and Seung \cite{Lee99} proposed and developed Nonnegative Matrix Factorization (NMF) algorithms, and demonstrated
that NMF has part-based representation which can be used for
intuitive perception interpretation.

\subsection{Related Work}

NMF has emerged in 1994 by Paatero and Tapper \cite{paatero1994positive}
for performing environmental data analysis.
The purpose of NMF is to decompose an input $m$-by-$n$ nonnegative matrix $A \in \mathbb{R}_{+}^{m\times n}$
into $m$-by-$r$ nonnegative matrix $B \in \mathbb{R}_{+}^{m\times r}$ and
$r$-by-$n$ nonnegative matrix $C \in \mathbb{R}_{+}^{r\times n}$:
$$
A \approx B C,
$$
and
\begin{equation} \label{nmfa}
\min_{ B, C \geq 0} \ \| A - B C \|^2_F,
\end{equation}
where $B, C \geq 0$ means that each entry of $B$ and $C$ is nonnegative,
$\| \cdot \|_F$ is the Frobenius norm of a matrix,
and $r$ (the low rank value) is smaller than $m$ and $n$. For simplicity, we assume that $m \ge n$,
Several researchers have proposed and developed
algorithms for determining such nonnegative matrix factorization in the literature.
For instance, Lee and Seung \cite{Lee99,Lee01} proposed to solve NMF by using
the multiplicative update algorithm by finding both $B$ and $C$ iteratively.
Also Yuan and Oja \cite{Yuan05} considered and studied
a projective nonnegative matrix factorization and proposed the following
minimization problem:
$$
\min_{ B \geq 0} \ \| A - B B^T A \|^2_F,
$$
where $B^T$ is the transpose of $B$.
In the optimization problem, it is required to
find a projection matrix $BB^T$ such that the difference between the given nonnegative matrix $A$ and
its projection $BB^T A$ is as small as possible.

We remark that there can be many possible solutions in (\ref{nmfa}). In practice,
it is necessary to impose additional constraints for finding NMF.
In some applications, orthogonality, sparsity and/or smoothness constraints on $B$ and/or $C$
are incorporated in (\ref{nmfa}). Because of these constraint formulations, many optimization techniques
have been designed to solve these minimization problems. For example,
the multiplicative update algorithms \cite{Choi08,Cich09,Ding06} are revised
to deal with  these constraint minimization problems.

\subsection{The Contribution}

The proposed Nonnegative Low Rank Matrix (NLRM) approximation method is
completely different from classical NMF.
Here NLRM approximation is to find a nonnegative low rank matrix $X$ such that
$X \approx A$ such that their difference is as small as possible.
Mathematically, it can be formulated as the following optimization problem
\begin{equation}\label{pmain}
\min_{\rank({X})=r,{ X}\geq 0} \ \| { A}- { X}\|_{\textrm{F}}^{2}.
\end{equation}
There are two advantages in the proposed NLRM method.
\begin{itemize}
\item
It is obvious in (\ref{nmfa}) that
when $B$ and $C$ are nonnegative, then the resulting matrix $BC$ is also nonnegative.
But these constraints are more restricted than that required in (\ref{pmain}).
Instead of using NMF in (\ref{nmfa}), we study NLRM in (\ref{pmain}).
The distance $\| { A}- { X}\|_{\textrm{F}}^{2}$
by the proposed NLRM method can be smaller than
$\min_{ B, C \geq 0} \ \| A - B C \|^2_F$ by the NMF method.
It implies that the proposed NLRM method can obtain a better low rank matrix
approximation.
\item
The proposed NLRM approximation
admits a matrix singular value decomposition, i.e.,
\begin{equation} \label{svd}
X = U \Sigma V^T,
\end{equation}
where $U$ is an $m$-by-$n$ matrix,
$\Sigma$ is an $n$-by-$n$ diagonal matrix,
and $V^T$ is also an $n$-by-$n$ matrix.
The columns of $U$
are called the left singular vectors of the singular value decomposition $\{ u_i \}_{i=1}^{m}$.
These left singular vectors form an orthonormal basis system in $\mathbb{R}^{m\times m}$
such that $u_i^T u_j = 1$ if $i = j$, otherwise 0.
The rows of $V^T$ refer to the elements of the right singular vectors of the singular value decomposition
$\{ v_i \}_{i=1}^{n}$. These right singular vectors also form an orthonormal basis system in
$\mathbb{R}^{n\times n}$ such that $v_i^T v_j = 1$ if $i = j$, otherwise 0.
The diagonal elements of $\Sigma = \ {\rm diag}( \sigma_1, \sigma_2, \cdots, \sigma_n)$
are called the singular values. As $X$ is a rank $r$ matrix, we have
$\sigma_i \ge 0$ for $1 \le i \le r$ and $\sigma_i = 0$ for $r+1 \le i \le n$.
The ordering of the singular values follows the descending order, i.e.,
$\sigma_1 \ge \sigma_2 \ge \cdots, \ge \sigma_r$.
We remark that both $U$ and $V$ are not necessary to be nonnegative, but the resulting
matrix $X = U \Sigma V^T$ must be non-negative.
Therefore, this decomposition is different from principal component analysis as there is such requirement in
principal component analysis. According to the singular value decomposition of $X$,
the proposed method can
identify important singular basis vectors based on singular values.
In the classical NMF, this information
cannot be obtained directly.
\end{itemize}
Our experimental results are presented to demonstrate the above mentioned advantages of the proposed NLRM method
compared the NMF method.

The paper is organized as follows. In Section 2, we present our algorithm and show the convergence.
In Section 3, numerical results are presented to demonstrate the proposed algorithm.
Finally, some concluding remarks and future research work are given in Section 4.

\section{The Optimization on Manifolds}

Constrained optimization is quite well studied as an area of research, and many powerful methods are proposed to solve the general problems in that area. In some special cases, the constrain sets possess particularly rich geometric properties, i.e., they are manifolds in the meaning of classical differential geometry. Then some constrained optimization problems can be rewritten as optimizing a real-valued function $f(x)$ on a manifold ${\cal{M}}$:
\begin{align}\label{opm}
\min_{x\in {\cal M}} f(x).
\end{align}
Here, ${\cal M}$ can be the Stiefel manifold, the Grassmann manifold and the fixed rank matrices manifold and so on.  In order to better understand  manifolds and some related definitions, e.g., charts, atlases and tangent spaces, we refer to \cite{lee2013smooth} and the references therein.  In general, the dimensions that some classical constrained techniques work are much bigger than the corresponding manifold (see e.g., \cite{absil2009optimization,vandereycken2013low}).

\subsection{The Algorithm}

Alternating projection method is popular in searching a point in the intersection of convex sets because of its simplicity and intuitive appeal. Its basic idea is  iteratively projecting a point one set and then the other. In contrast to the well known cases, this paper concerns with the extensions of convex sets to non-convex sets. Here,  one set is the $m\times n$ fixed rank matrices manifold
\begin{align}\label{v1}
\mathcal{M}_{r}:=\left\{ X\in \mathbb{R}^{m\times n}, \rank(X)= r\right\},
\end{align}
and the other one is the convex set of $m\times n$ nonnegative matrices
\begin{align}
\mathcal{M}_{n}:=\left\{X\in \mathbb{R}^{m\times n}, X_{i,j}\geq 0, i=1,...,m,j=1,...,n\right\}.\label{v2}
\end{align}
In order to introduce the main algorithm, we need to define two projections that project the given matrix onto $\mathcal{M}_{r}$ and $\mathcal{M}_{n}$, respectively. By the Eckart-Young-Mirsky theorem \cite{golub2012matrix}, the projection onto fixed rank matrix set $\mathcal{M}_{r}$ can be expressed as follows:
\begin{align}\label{p1}
\pi_{1}({ X})=\sum_{i=1}^{r}\sigma_{i}(X) u_{i}(X) {v}_{i}^{T}(X),
\end{align}
where $\sigma_{i}(X)$ is the $i$-th singular values of $X$, and their corresponding left and right singular vectors:
$u_{i}(X)$ and $v_{i}(X)$.
The projection onto the nonnegative matrix set $\mathcal{M}_{n}$ can be expressed as follows:
\begin{align}\label{p2}
\pi_{2}({ X})=\left\{\begin{array}{c}
                     X_{ij}, ~~~ {\rm if} ~~ X_{ij}\geq 0, \\
                     0,     ~~~~~ {\rm if} ~~   X_{ij} < 0.
                   \end{array}
\right.
\end{align}
In particular, since the manifold $\mathcal{M}_{r}$ is not convex,
the projection mapping $\pi_{1}(X)$ can no longer be single valued (for example when the $r$-th singular value has multiplicity higher than 1).
Then the algorithm of alternating projection on $\mathcal{M}_{r}$ and $\mathcal{M}_{n}$ can be given as Algorithm \ref{ag1}. The framework of this algorithm is the same as the general case, while the only difference is that the projections are respectively chosen as $\pi_{1}$ and $\pi_{2}$ given in \eqref{p1} and \eqref{p2}.

\begin{algorithm}[h]
\caption{Alternating Projections On Manifolds} \label{ag1}
\textbf{Input: } Given a nonnegative matrix ${A}\in \mathbb{R}^{m\times n}$ this algorithm computes nearest rank-$r$ nonnegative matrix. \\
~~1: Initialize ${ X}_0 =  A$; \\
~~2: \textbf{for} $k=1,2,...$\\
~~3: \quad ${ Y}_{k+1}=\pi_{1}({ X}_{k});$\\
~~4: \quad ${ X}_{k+1}=\pi_{2}({ Y}_{k+1});$\\
~~5: \textbf{end}\\
\textbf{Output:} ${ X}_k$ when the stopping criterion is satisfied.
\end{algorithm}

We note in Algorithm 1 that the projection of a given matrix onto the manifold $\mathcal{M}_r$ is done by truncating
small singular values in the singular value decomposition of the given matrix. The computational complexity
of this procedure is of $(mnr)$ operations.

Different from the convex sets case, the intersection of $\mathcal{M}_{r}$ and $\mathcal{M}_{n}$ decides whether the sequence generated by Algorithm \ref{ag1}, converges or not.
In order to show the convergence of Algorithm \ref{ag1}, we
compute the dimension of the intersection of $\mathcal{M}_{r}$ and $\mathcal{M}_{n}$.
%In order to show the convergence of Algorithm \ref{ag1}, we need to establish the following results.

\begin{theorem}\label{th1}
Let $\mathcal{M}_{r}$ and $\mathcal{M}_{n}$ be defined as \eqref{v1}-\eqref{v2}. Then
\begin{equation}\label{v3}
\mathcal{M}_{rn}:=\mathcal{M}_{r}\cap \mathcal{M}_{n}=\left\{X\in \mathbb{R}^{m\times n}, ~ \rank(X)=r,~X_{ij}\geq 0,~ i=1,...,m,j=1,...,n\right\}
\end{equation}
is a smooth manifold with dimension  $(m+n)r-r^2$.
\end{theorem}

The proof of Theorem \ref{th1} can be found in Supplementary.
Moreover, we need to define the angle $\alpha(A)$ of $A\in \mathcal{M}_{rn}$
where
$$
\alpha(A)=cos^{-1}(\sigma(A)) \quad {\rm and} \quad
\sigma(A)=\lim_{\xi\rightarrow 0} \sup_{B_{1}\in F^{\xi}_{1}(A), B_{2}\in F^{\xi}_{2}(A)}
\left\{\frac{\left<B_{1}-A,B_{2}-A\right>}{\|B_{1}-A\|_{F}\|B_{2}-A\|_{F}}\right\},
$$
with
$$
F_{1}^{\xi}(A)=
\{
B_1 \ | \ B_1 \in \mathcal{M}_{r}\backslash A,
\|B_1-A\|_{F}\leq \xi, B_{1}-A \bot T_{\mathcal{M}_{r}\cap \mathcal{M}_{n}}(A)
\},
$$
$$
F_{2}^{\xi}(A)=
\{
B_2 \ | \ B_2 \in \mathcal{M}_{n}\backslash A,
\|B_2-A\|_{F}\leq \xi, B_{2}-A \bot T_{\mathcal{M}_{r}\cap \mathcal{M}_{n}}(A)
\},
$$
and
$T_{\mathcal{M}_{r}\cap \mathcal{M}_{n}}(A)$ is
the tangent space of
$\mathcal{M}_{r} \cap \mathcal{M}_{n}$ at point $A$ (the definition of the tangent space can be found in Supplemary).
The angle is calculated based on the two points
belonging $\mathcal{M}_1$ and $\mathcal{M}_2$ respectively.

A point $A$ in $\mathcal{M}_{rn}$
is nontangential if $\alpha(A)$ has a positive angle \cite{andersson2013alternating}, i.e.,
$0 \le \sigma(A) < 1$. It is interesting to note that if there is one
nontangential point, majority of points are also nontangential because of the manifold smoothness.
We show that $\mathcal{M}_{rn}$ contains a non-empty set $\mathcal{M}^{nt}_{rn}$
of nontangential points. The proof can be found in Supplementary.

\begin{theorem}\label{theorem2}
$\mathcal{M}^{nt}_{rn} \neq \emptyset$.
\end{theorem}

By using this property, Andersson and Carlsson \cite{andersson2013alternating} have shown the following theorem.

\begin{theorem}\label{thm_convergence}
Let $\mathcal{M}_{r}$, $\mathcal{M}_{n}$ and $\mathcal{M}_{rn}$ be given as \eqref{v1}, \eqref{v2} and \eqref{v3}, the projections onto $\mathcal{M}_{r}$ and $\mathcal{M}_{rn}$ be given as \eqref{p1}-\eqref{p2}, respectively. Denote
$\pi$ as the projection onto $\mathcal{M}_{rn}$.
Suppose that $P \in \mathcal{M}^{nt}_{rn}$ is a non-tangential intersection point, then for any  given $\epsilon>0$ and
$1>c>\sigma(P)$, there exist an
$\xi>0$ such that for any ${ A}\in Ball({P},\xi)$ (the ball neighborhood of $P$ with radius $\xi$ contains the given
nonnegative matrix $A$)
  the sequence $\{{X}_k\}_{k=0}^\infty$ generated by the alternating projection algorithm initializing from given $A$ satisfies the following results:
  \begin{enumerate}[(1)]
\item the sequence converges to a point $X_\infty \in {\cal M}_{rn}$,
\item $\| { X}_\infty - \pi({ A}) \|_{F} \leq \epsilon \| {A} - \pi({ A}) \|_{F}$,
\item $\| { X}_\infty - {X}_k \|_{F} \leq {\rm const} \cdot c^k \| {A} - \pi({ A}) \|_{F}$,
  \end{enumerate}
where $\pi(A) \equiv {\rm argmin}_{\rank({X})=r,{ X}\geq 0} \ \| { A}- { X}\|_{\textrm{F}}^{2}$ (the optimized solution).
\end{theorem}

According to Theorem \ref{thm_convergence}, we know Algorithm 1 can converge and find a matrix that can be sufficiently close
the best nonnegative approximation $\pi({ A}) $.

\section{Numerical Results}

In this section, we conducted several experiments to demonstrate that the
proposed NLRM method can obtain a better low rank matrix
approximation, and can provide
a significant index based on singular values that can be used to
identify important singular basis vectors in the approximation.

\subsection{The First Experiment}

In the first experiment, we randomly generated $m$-by-$r$ nonnegative matrices $B$ and
$r$-by-$n$ nonnegative matrices $C$ where their matrix entries follow
a uniform distribution in between 0 and 1. The given non-negative matrix $A$ was computed by
the multiplication of $B$ with $C$.
We employed the proposed NLRM algorithm (Algorithm 1) to test the relative residual $\| A - X_{c} \|_F / \| A \|_F$
and compared with testing NMF algorithms:
A-MU \cite{gillis2012accelerated}, HALS \cite{cichocki2007hierarchical}, A-HALS \cite{gillis2012accelerated}, PG \cite{Lin06a} and A-PG \cite{Lin06a}.
Here $X_c$ are the computed solutions by different algorithms.

Tables \ref{table1a}-\ref{table1c} shows the relative
residuals of the computed solutions from the proposed algorithm and the testing NMF algorithms
for synthetic data sets of sizes 100-by-80, 200-by-160 and 500-by-400 respectively.
When there is no noise (the second column in the tables), the input data matrix $A$ is exactly a nonnegative rank $r$ matrix.
Therefore, the proposed NLRM algorithm can provide exact recovery results in the first iteration.
However, there is no guarantee that
other testing NMF algorithms can determine the underlying nonnegative low rank factorization.
In the tables, it is clear that the testing NMF algorithms cannot obtain the underlying low rank factorization.
One of the reason may be that NMF algorithms can be sensitive to initial guesses.
In the tables, we illustrate this phenomena by displaying
the mean relative residual and the range containing both the minimum and the maximum relative residuals
by using ten initial guesses randomly generated.
We find in the tables that the relative residuals computed by the proposed
NLRM algorithm are always smaller than the minimum relative residuals by the testing NMF algorithms.

On the other hand, the performance of the proposed algorithm was evaluated when
a Gaussian noise of zero mean and variance $\sigma$ ($=0.001,0.005,0.01$) were added in the generated matrices $A$.
The relative residuals of the computed solutions by the proposed NLRM algorithm and the other NMF algorithms are reported
in Tables \ref{table1a}-\ref{table1c}. We see from the tables that the relative residuals computed by the proposed
method are always smaller than the minimum relative residuals by the testing NMF algorithms.
All these results show that the performance of the proposed NLRM algorithm is better than that of the testing NMF algorithms.

\begin{table}
 \begin{center}
  \setlength{\abovecaptionskip}{-1pt}
       \setlength{\belowcaptionskip}{-1pt}
\caption{The relative residuals $\| A - X_{c} \|_F / \| A \|_F$ by different algorithms for 100-by-80 synthetic data matrices with
several noise levels.}\label{table1a}
 \medskip
\scalebox{0.8}{
\begin{tabular}{|l||c|c|c|c|}
\hline
       & \multicolumn{4}{c|}{$r=10$} \\
\cline{2-5}
Method  & no noise & $\sigma=0.001$ & $\sigma=0.005$ & $\sigma=0.01$ \\
\hline
NLRM   &   9.15e-16        &  1.72e-04            &  8.62e-04            &  1.70e-03  \\
\hline
A-MU (mean)   &  1.43e-03         &  1.43e-03           &  1.68e-03          &  2.26e-03 \\
A-MU (range)  &  [8.66e-04, 2.30e-03]  &  [8.85e-04, 2.30e-03] & [1.20e-03, 2.50e-03] & [1.90e-03, 2.90e-03]  \\
\hline
HALS (mean)   &  1.36e-03        &  1.38e-03     &  1.63e-03            &  2.22e-03 \\
HALS (range)  &  [7.22e-04, 2.30e-03] & [7.33e-04, 2.30e-03]   & [1.10e-03, 2.40e-03] & [1.80e-03, 2.90e-03]  \\
\hline
A-HALS (mean)   &  1.25e-03        &       1.37e-03      &    1.61e-03        &  2.23e-03 \\
A-HALS (range)  & [6.84e-04, 2.40e-03]     &      [5.22e-04, 2.30e-03]        &     [5.22e-04, 2.30e-03]         &  [6.84e-04, 2.40e-03] \\		
\hline
PG (mean)   &    	1.06e-03 &      8.24e-04      &      1.26e-03        &  1.96e-03 \\
PG (range)  &   [2.36e-04, 2.00e-03]         &    [3.03e-04, 2.00E-03]	  &    [9.39e-04, 1.80e-03]	&  [1.80e-03, 2.20e-03]\\	
\hline
A-PG (mean)   &    	2.03e-03        &  2.04e-03          & 2.21e-03           & 2.66e-03 \\
A-PG (range)  &    [1.40e-03, 2.90e-03]       &    [1.40e-03, 2.90e-03]          &  [1.60e-03, 3.00e-03]          &  [2.20e-03, 3.40e-03]	 \\		
\hline
\hline
       & \multicolumn{4}{c|}{$r=20$} \\
\cline{2-5}
Method  & no noise & $\sigma=0.001$ & $\sigma=0.005$ & $\sigma=0.01$ \\
\hline
NLRM   &    8.32e-16        & 1.22e-04                     &  6.07e-04                    &   1.20e-03         \\
\hline
A-MU (mean)   & 2.20e-04                         & 2.54e-04                   & 6.49e-04                 & 1.21e-03         \\
A-MU (range)  & [9.74e-05, 4.41e-04]             & [1.56e-04, 4.51e-04]       & [6.15e-04, 7.36e-04]     & [1.20e-03, 1.30e-03]    \\		
\hline
HALS (mean)   &2.16e-04   &   2.53e-04                  &   6.51e-04                   &   1.22e-03        \\
HALS (range)  &[5.93e-05, 5.00e-04]   &  [1.36e-04, 4.98e-04]   & [6.10e-04, 7.47e-04] & [1.20e-03, 1.30e-03]   \\
\hline
A-HALS (mean)   &  5.47e-05  & 1.42e-04    &6.12e-04    &1.20e-03 \\
A-HALS (range)  & [2.29e-05, 1.19e-04] &[1.25e-04, 1.88e-04]     &[6.09e-04, 6.20e-04]  &[1.20e-03, 1.20e-03]\\
\hline
PG (mean)   & 1.15e-04    & 1.40e-04     & 6.17e-04    &1.20e-03 \\
PG (range)  & [4.33e-05, 1.76e-04]   &[1.22e-04, 1.90e-04]       &[6.08e-04, 6.41e-04]     &[1.20e-03, 1.20e-03] \\
\hline
A-PG (mean)   &  5.23e-04   &  5.37e-04     &8.06e-04     & 1.32e-03  \\
A-PG (range)  & [3.30e-04, 6.98e-04]   & [3.51e-04, 7.10e-04]     & [6.90e-04, 9.30e-04]  &[1.30e-03, 1.40e-03]\\
\hline
\hline
 & \multicolumn{4}{c|}{$r=40$} \\
\cline{2-5}
Method  & no noise & $\sigma=0.001$ & $\sigma=0.005$ & $\sigma=0.01$ \\
\hline
NLRM   &  2.61e-15   &   8.43e-05      &    4.21e-04     &  8.43e-04          \\
\hline
A-MU (mean)   &   1.82e-03     &  1.84e-03        & 1.92e-03      & 2.03e-03\\
A-MU (range)  &   [1.20e-03, 2.30e-03]         &  [1.20e-03, 2.30e-03]        & [1.30e-03, 2.50e-03]      &  [1.60e-03, 2.40e-03]        \\
\hline
HALS (mean)   &2.82e-03 &  2.83e-03        &  2.87e-03      &3.00e-03  \\
HALS (range)  & [2.40e-03, 3.20e-03]&  [2.40e-03, 3.30e-03]    & [2.50e-03, 3.30e-03]   & [2.60e-03, 3.40e-03]\\
\hline
A-HALS (mean)   &  1.33e-05    & 8.64e-05      & 4.23e-04      &8.43e-04   \\
A-HALS (range)  &  [1.31e-06, 7.48e-05]   &  [8.43e-05, 9.35e-05]    &[4.21e-04, 4.35e-04]     & [8.43e-04, 8.44e-04]  \\
\hline
PG (mean)   &4.39e-04   &  2.42e-04    & 4.84e-04  &1.02e-03 \\
PG (range)  & [4.07e-05, 7.64e-04]    &[9.74e-05, 7.83e-04]      &  [4.23e-04, 6.98e-04]    & [8.44e-04, 1.70e-03]\\
\hline
A-PG (mean)   & 3.40e-04    &  3.64e-04    & 5.87e-04  &9.48e-04\\
A-PG (range)  &  [2.15e-06, 1.00e-03]   &   [8.43e-05, 1.00e-03]   &[4.21e-04, 1.10e-03]      &[8.43e-04, 1.30e-03] \\
\hline
\end{tabular}}
\end{center}
\end{table}

\begin{table}
\begin{center}
 \setlength{\abovecaptionskip}{-1pt}
       \setlength{\belowcaptionskip}{-1pt}
\caption{The relative residuals $\| A - X_{c} \|_F / \| A \|_F$ by different algorithms for 200-by-160 synthetic data matrices with
several noise levels.}\label{table1b}
\medskip
\scalebox{0.8}{
\begin{tabular}{|l||c|c|c|c|}
\hline
       & \multicolumn{4}{c|}{$r=10$} \\
\cline{2-5}
Method  & no noise & $\sigma=0.001$ & $\sigma=0.005$ & $\sigma=0.01$ \\
\hline
NLRM   &   5.25e-16        &  1.33e-04          &  6.67e-04            &  1.30e-03  \\
\hline
A-MU (mean)   &  3.38e-03         &  3.39e-03           &   3.45e-03           & 3.65e-03 \\
A-MU (range)  &  [1.70e-03, 4.90e-03]  &  [1.80e-03, 4.90e-03] & [1.90e-03, 4.90e-03] & [2.30e-03, 5.00e-03]  \\
\hline
HALS (mean)   &  4.49e-03          &  4.48e-03    &  4.47e-03            & 4.61e-03 \\
HALS (range)  &  [3.30e-03, 6.00e-03] & [3.30e-03, 5.90e-03]   & [3.40e-03, 5.60e-03] & [3.50e-03, 5.70e-03]  \\
\hline
A-HALS (mean)   &  4.07e-03          & 4.23e-03     & 4.10e-03   &4.05e-03 \\
A-HALS (range)  &  [2.30e-03, 5.70e-03]   & [3.00e-03, 5.30e-03]  & [2.80e-03, 5.60e-03]     & [2.50e-03, 5.70e-03]\\
\hline
PG (mean)   &  2.79e-03   &  2.09e-03  & 1.97e-03   &2.60e-03\\
PG (range)  & [8.70e-04, 4.10e-03]  &  [2.67e-04, 4.10e-03]    & [1.00e-03, 3.40e-03]    &[1.50e-03, 4.40e-03]   \\
\hline
A-PG (mean)   & 3.02e-03   & 3.03e-03  &  3.24e-03  &3.39e-03\\
A-PG (range)  & [6.98e-04,	5.20e-03]  & [7.09e-04, 5.20e-03]   &[9.58e-04, 5.20e-03] & [1.50e-03, 5.40e-03] \\
\hline
\hline
       & \multicolumn{4}{c|}{$r=20$} \\
\cline{2-5}
Method  & no noise & $\sigma=0.001$ & $\sigma=0.005$ & $\sigma=0.01$ \\
\hline
NLRM   & 6.41e-15   &   9.24e-05      & 4.62e-04  & 9.25e-04  \\
\hline
A-MU (mean)   &   1.73e-03   &  1.73e-03      &   1.80e-03       & 1.97e-03   \\
A-MU (range)  & [1.40e-03, 2.50e-03]  &  [1.40e-03, 2.50e-03]   & [1.50e-03, 2.60e-03]     &  [1.70e-03, 2.70e-03]   \\
\hline
HALS (mean)   & 1.55e-03  &     1.55e-03      &  1.63e-03       & 1.82e-03    \\
HALS (range)  & [1.40e-03, 1.70e-03]& [1.40e-03, 1.70e-03]   &   [1.50e-03, 1.80e-03]     &  [1.70e-03, 2.00e-03]  \\
\hline
A-HALS (mean)   &  1.29e-03 &  1.39e-03  & 1.40Ee-03    & 1.63e-03 \\
A-HALS (range)  &  [1.00e-03, 1.70e-03] &  [1.10e-03, 1.80e-03]   & [8.55e-04, 1.80e-03]  & [1.40e-03, 1.90e-03]  \\
\hline
PG (mean)   &  1.80e-03   & 1.92e-03   &  1.73e-03    &1.88e-03   \\
PG (range)  &  [1.40e-03, 2.30e-03]  &[1.60e-03, 2.10e-03]   &  [1.50e-03, 2.10e-03]  & [1.60e-03, 2.30e-03]  \\
\hline
A-PG (mean)   & 9.84e-04 & 9.86e-04  &  1.10e-03  & 1.35e-03 \\
A-PG (range)  &  [7.59e-04, 1.20e-03]  & [7.65e-04, 1.20e-03] & [9.02e-04, 1.30e-03]   & [1.20e-03, 1.50e-03]  \\
\hline
\hline
       & \multicolumn{4}{c|}{$r=40$} \\
\cline{2-5}
Method  & no noise & $\sigma=0.001$ & $\sigma=0.005$ & $\sigma=0.01$ \\
\hline
NLRM   &  8.76e-16    &  6.28e-05   &  3.14e-04&  6.28e-04  \\
\hline
A-MU (mean)   &  6.05e-04    &  6.09e-04    &  6.84e-04   &  8.74e-04   \\
A-MU (range)  &  [5.26e-04, 6.61e-04]    &   [5.31e-04, 6.64e-04]   & [6.16e-04, 7.30e-04]   & [8.24e-04, 9.08e-04]    \\
\hline
HALS (mean)   & 3.54e-03  &  3.56e-03   & 3.63e-03     &  3.72e-03   \\
HALS (range)  & [2.90e-03, 4.80e-03] &  [2.90e-03, 4.80e-03]    &  [2.90e-03, 4.90e-03]     & [2.90e-03, 5.10e-03]   \\
\hline
A-HALS (mean)   &  2.37e-05    & 6.88e-05  & 3.15e-04  &  6.28e-04 \\
A-HALS (range)  &  [1.51e-05, 4.15e-05]    & [6.37e-05, 8.13e-05]   &  [3.14e-04, 3.17e-04]    & [6.28e-04, 6.29e-04] \\
\hline
PG (mean)   &  3.05e-03   &  2.84e-03  &3.21e-03    &3.22e-03   \\
PG (range)  & [1.70e-03,4.70e-03]  & [1.30e-03, 6.10e-03] & [1.40e-03, 5.70e-03]   & [1.50e-03, 5.20e-03]  \\
\hline
A-PG (mean)   &  6.05e-04  &  6.09e-04   &  6.83e-04 & 8.74e-04 \\
A-PG (range)  &  [5.04e-04, 7.43e-04]    & [5.07e-04, 7.46e-04]    &[5.93e-04, 8.06e-04]    &[8.05e-04, 9.71e-04]   \\
\hline
\end{tabular}}
\end{center}
\end{table}

\begin{table}
\begin{center}
\setlength{\abovecaptionskip}{-1pt}
       \setlength{\belowcaptionskip}{-1pt}
\caption{The relative residuals $\| A - X_{c} \|_F / \| A \|_F$ by different algorithms for 500-by-400 synthetic data matrices
with several noise levels.}\label{table1c}
\medskip
\scalebox{0.8}{
\begin{tabular}{|l||c|c|c|c|}
\hline
       & \multicolumn{4}{c|}{$r=10$} \\
\cline{2-5}
Method  & no noise & $\sigma=0.001$ & $\sigma=0.005$ & $\sigma=0.01$ \\
\hline
NLRM   &   1.33e-15         &  8.07e-05           &  4.03e-04           &  8.07e-04 \\
\hline
A-MU (mean)   &  1.21e-03          &  1.22e-03           &   1.28e-03         &  1.47e-03 \\
A-MU (range)  &  [6.06e-04, 3.00e-03]  &  [6.14e-04, 3.00e-03] & [7.38e-04, 3.10e-03] & [1.00e-03, 3.20e-03]  \\
\hline
HALS (mean)   &  4.05e-03          &  4.09e-03     &  4.14e-03           & 4.34e-03 \\
HALS (range)  &  [1.00e-03, 6.60e-03] & [1.10e-03, 6.60e-03]   & [1.10e-03, 6.70e-03] & [1.30e-03, 8.40e-03]  \\
\hline
A-HALS (mean)   &   3.39e-03   & 4.33e-03   &  3.45e-03   &4.14e-03   \\
A-HALS (range)  & [1.10e-03, 5.80e-03]   &[8.69e-04, 6.60e-03]  & [8.58e-04, 5.00e-03]   &[1.20e-03, 6.90e-03] \\
\hline
PG (mean)   &    5.17e-03   & 5.38e-03   &  4.52e-03  & 5.71e-03  \\
PG (range)  & [2.60e-03, 7.60e-03]  &[2.40e-03, 7.40e-03] &[1.50e-03, 6.20e-03]   & [2.40e-03, 7.30e-03]  \\
\hline
A-PG (mean)   &   3.26e-04 & 3.74e-04  & 4.22e-04    & 8.15e-04  \\
A-PG (range)  & [3.93e-05, 2.30e-03]   & [9.01e-05, 2.50e-03]    & [4.06e-04, 4.54e-04]   &[8.08e-04, 8.31e-04]   \\
\hline
\hline
       & \multicolumn{4}{c|}{$r=20$} \\
\cline{2-5}
Method  & no noise & $\sigma=0.001$ & $\sigma=0.005$ & $\sigma=0.01$ \\
\hline
NLRM   &  1.56e-15   & 5.81e-05   &  2.91e-04  &  5.81e-04   \\
\hline
A-MU (mean)   &4.38e-03  &4.38e-03      &  4.38e-03     &  4.42e-03   \\
A-MU (range)  &  [4.10e-03, 4.80e-03]     & [4.10e-03, 4.80e-03]    & [4.10e-03, 4.80e-03]	    & [4.20e-03, 4.80e-03]         \\
\hline
HALS (mean)   &4.67e-03   &  4.67e-03       &   4.69e-03     &  4.73e-03  \\
HALS (range)  & [4.40e-03, 5.00e-03]  &  [4.40e-03, 5.00e-03]    &   [4.40e-03, 5.00e-03]    & [4.40e-03, 5.10e-03]	  \\
\hline
A-HALS (mean)   &  4.32e-03   & 4.29e-03   &4.48e-03   &4.41e-03   \\
A-HALS (range)  &  [3.80e-03, 4.70e-03]  & [4.10e-03, 4.50e-03]    & [4.20e-03, 4.90e-03]   &[4.10e-03, 4.80e-03]	   \\
\hline
PG (mean)   &   8.87e-03  &  9.03e-03 &   8.82e-03   &  8.82e-03 \\
PG (range)  & [8.10e-03, 9.60e-03]   &[8.30e-03, 1.01e-02]   & [8.10e-03, 9.60e-03]  &[8.10e-03, 9.40e-03]   \\
\hline
A-PG (mean)   &    4.58e-03  &  4.58e-03    &   4.59e-03  & 4.61e-03\\
A-PG (range)  &  [4.10e-03, 5.10e-03]   &  [4.10e-03, 5.10e-03]   & [4.10e-03, 5.10e-03]  & [4.20e-03, 5.10e-03]\\
\hline
\hline
       & \multicolumn{4}{c|}{$r=40$} \\
\cline{2-5}
Method  & no noise & $\sigma=0.001$ & $\sigma=0.005$ & $\sigma=0.01$ \\
\hline
NLRM   &  1.78e-15   &   4.10e-05     & 2.05e-04    &4.10e-04  \\
\hline
A-MU (mean)   &   2.10e-03    & 2.10e-03    &  2.10e-03        &   2.11e-03  \\
A-MU (range)  &   [1.90e-03, 2.30e-03]     &   [1.90e-03, 2.30e-03]  &   [1.90e-03, 2.30e-03]  & [1.90e-03, 2.30e-03]   \\
\hline
HALS (mean)   &5.58e-03 &5.58e-03    &    5.51e-03     & 5.46e-03    \\
HALS (range)  &  [4.50e-03, 7.10e-03] &[4.50e-03, 7.10e-03]      &  [4.50e-03, 7.20e-03]   &   [4.20e-03, 7.20e-03]   \\
\hline
A-HALS (mean)   &  2.05e-03   &  1.98e-03  &  2.26e-03     & 1.78e-03  \\
A-HALS (range)  & [1.60e-03, 3.30e-03] & [1.40e-03, 2.70e-03]   &  [1.50e-03, 5.00e-03]   & [1.50e-03, 2.60e-03]	  \\
\hline
PG (mean)   & 7.53e-03 &  7.83e-03    &  7.74e-03& 7.77e-03  \\
PG (range)  &  [6.60e-03, 8.90e-03]   & [6.80e-03, 8.90e-03]  &  [6.70e-03, 9.00e-03] & [6.60e-03, 9.70e-03]  \\
\hline
A-PG (mean)   &  2.55e-03   &    2.55e-03          & 2.56e-03      & 2.60e-03  \\
A-PG (range)  &  [2.50e-03, 2.70e-03]  & [2.50e-03, 2.70e-03]  & [2.50e-03, 2.70e-03]  &[2.50e-03, 2.70e-03]   \\
\hline
\end{tabular}}
\end{center}
\end{table}

\subsection{The Second Experiment}

In the second experiment, we randomly generated $m$-by-$n$ nonnegative matrices $A$
where their matrix entries follow
a uniform distribution in between 0 and 1. The low rank minimizer is unknown
in this setting. Table \ref{table1d} shows that the relative residuals $\| A - X_{c} \|_F / \| A \|_F$
of the computed solution $X_{c}$ from the proposed NLRM algorithm and the testing NMF algorithms.
In the testing NMF algorithms, we use 10 different initial guesses and report the results of
the mean and the range (the minimum and the maximum values) of the relative residuals in the tables.
We see from Table \ref{table1d} that the relative residuals computed by the proposed
NLRM method are smaller than the minimum relative residuals by the testing NMF algorithms.
All these results show that the performance of the proposed NLRM algorithm is better than that of the testing NMF algorithms.

Moreover, we considered the CBCL face database \cite{face}.
In the face database, there are $m=2469$ facial images, each consisting of $n=19\times19 = 361$ pixels, and constituting a face image matrix $A\in \mathbb{R}_{+}^{361 \times 2469}$. We tested several values of $r=20,40,60,80$ for nonnegative low rank minimization
and compared the proposed NLRM algorithm with the testing NMF algorithms.
In the testing NMF algorithms, we used 10 different initial guesses and report the results of
relative residuals of the mean and the range in the tables.
We see from Table \ref{table1e} that the relative residuals computed by the proposed
NLRM method are smaller than the minimum relative residuals by the testing NMF algorithms.
Again the performance of the proposed algorithm is better than that of the testing NMF algorithms.

\begin{table}
\begin{center}
\setlength{\abovecaptionskip}{-1pt}
       \setlength{\belowcaptionskip}{-1pt}
\caption{The relative residuals $\| A - X_{c} \|_F / \| A \|_F$ by different algorithms for synthetic data matrices.}\label{table1d}
\medskip
\scalebox{0.8}{
\begin{tabular}{|l||c|c|c|}
\hline
       & \multicolumn{3}{c|}{$m=100$, $n=80$} \\
\cline{2-4}
Method  & $r=10$ & $r=20$ & $r=40$ \\
\hline
NLRM       &  0.4047           &  0.3228          & 0.1874 \\
\hline
A-MU (mean)   &  0.4087          &   0.3434         &  0.2487 \\
A-MU (range)   &  [0.4086, 0.4089] & [0.3428, 0.3447] & [0.2477, 0.2502]  \\
\hline
HALS (mean)       &  0.4087     &  0.3426           & 0.2448 \\
HALS (range)   & [0.4086, 0.4088]   & [0.3421, 0.3431] & [0.2438, 0.2458]  \\
\hline
A-HALS (mean)     & 0.4087   &  0.3424   &0.2449   \\
A-HALS (range)  & [0.4086, 0.409]  & [0.342, 0.3427]   &[0.2437, 0.2462] \\
\hline
PG (mean)   & 0.4087   & 0.3422  & 0.2445  \\
PG (range)  &[0.4086, 0.4089] &[0.3418, 0.3426]   & [0.2425, 0.2465]  \\
\hline
A-PG (mean)   & 0.4088 & 0.3428  & 0.2450  \\
A-PG (range)   & [0.4086, 0.4091]    & [0.3420, 0.3436]   &[0.2442, 0.2463]   \\
\hline
\hline
       & \multicolumn{3}{c|}{$m=200$, $n=160$} \\
\cline{2-4}
Method  & $r=10$ & $r=20$ & $r=40$ \\
\hline
NLRM    & 0.4503   &  0.4054 & 0.3252   \\
\hline
A-MU (mean)   &0.4522      &0.4151     & 0.3586   \\
A-MU (range)   & [0.4521, 0.4523]    & [0.4148, 0.4156]	    & [0.358, 0.3591]         \\
\hline
HALS (mean)    &  0.4522      &  0.4149     & 0.3569  \\
HALS (range)   &  [0.4521, 0.4523]    &   [0.4147, 0.4151]    & [0.3565, 0.3574]	  \\
\hline
A-HALS (mean)    & 0.4522  &0.4149   &0.3568   \\
A-HALS (range)    & [0.4521, 0.4523]    & [0.4145, 0.4152]   &[0.3565, 0.3572]	   \\
\hline
PG (mean)   & 0.4521 &  0.4147   &0.3568\\
PG (range)    &[0.4521, 0.4522]   & [0.4145, 0.415]  &[0.3564, 0.3572]   \\
\hline
A-PG (mean)    & 0.45215   & 0.41478  & 	0.35741\\
A-PG (range)    &  [0.4521, 0.4522]   & [0.4146, 0.4151]  & [0.3569, 0.3582]\\
\hline
\hline
       & \multicolumn{3}{c|}{$m=500$, $n=400$} \\
\cline{2-4}
Method  & $r=10$ & $r=20$ & $r=40$ \\
\hline

NLRM   &   0.4803     & 0.4607   &0.4238  \\
\hline
A-MU (mean)   & 0.4807    &  0.4635        &   0.4357  \\
A-MU (range)    &   [0.4807, 0.4807]  &   [0.4634, 0.4635	]  & [0.4357, 0.4358]   \\
\hline
HALS (mean)  &0.4807   &    0.4633     & 0.4352    \\
HALS (range) &[0.4807, 0.4808]      &  [0.4633, 0.4634]   &   [0.4352, 0.4353]   \\
\hline
A-HALS (mean)    & 0.4807  &  0.4634     & 0.4352  \\
A-HALS (range) & [0.4807, 0.4807]   &  [0.4633, 0.4634]   & [0.4351, 0.4353]	  \\
\hline
PG (mean)    &  0.4807   & 0.4633& 0.4353 \\
PG (range)    & [0.4807, 0.4807]  &  [0.4633, 0.4634] & [0.4351, 0.4354]  \\
\hline
A-PG (mean)    &    0.4807          & 0.4634      & 0.4353  \\
A-PG (range)  & [0.4807	0.4807]  & [0.4633, 0.4635]  &[0.4352, 0.4354]   \\
\hline
\end{tabular}}
\end{center}
\end{table}

\begin{table}
\begin{center}
\setlength{\abovecaptionskip}{-1pt}
       \setlength{\belowcaptionskip}{-1pt}
\caption{The relative residuals $\| A - X_{c} \|_F / \| A \|_F$ by different algorithms for face data.}\label{table1e}
\medskip
\scalebox{0.8}{
\begin{tabular}{|l||c|c|c|c|}
\hline
       & \multicolumn{4}{c|}{$m=361$, $n=2469$} \\
\cline{2-5}
Method  & $r=20$ & $r=40$ & $r=60$ & $r=80$ \\
\hline
NLRM   &   0.1170       &  0.0839          &  0.0654         & 0.0529 \\
\hline
A-MU (mean)   &  0.1220         & 0.0911          &  0.0751         &  0.0639 \\
A-MU (range)  &  [0.1218, 0.1222]  &  [0.0907, 0.0914] & [0.0746, 0.0755] & [0.0632, 0.0651]  \\
\hline
HALS (mean)   &  0.1220         &  0.0902    &  0.0721          & 0.0598 \\
HALS (range)  &  [0.1218, 0.1224] & [0.09, 0.0905]   & [0.0719, 0.0723] & [0.0596, 0.0601]  \\
\hline
A-HALS (mean)   &   0.1220   & 0.0901  & 0.0719   &0.0595   \\
A-HALS (range)  & [0.1217, 0.1222]   &[0.0898, 0.0904]  & [0.0717, 0.072]   &[0.0594, 0.0597] \\
\hline
PG (mean)   &    0.1219   & 0.0902   & 0.0740  & 0.0633 \\
PG (range)  & [0.1217, 0.122]  &[0.0899, 0.0909] &[0.0731, 0.0747]   & [0.0626, 0.0645]  \\
\hline
A-PG (mean)   &   0.1221 & 0.0908 &0.0754  & 0.0658  \\
A-PG (range)  & [0.1218, 0.1226]   & [0.0905, 0.0913]    & [0.075, 0.0757]   &[0.0652, 0.0684]   \\
\hline
\end{tabular}}
\end{center}
\end{table}

\subsection{The Third Experiment}

In practice, it is necessary to determine the value of rank for nonnegative matrix approximation.
In this experiment, we show the advantage of the proposed NLRM algorithm
for providing a significant index based on singular values that can be used to
identify important singular basis vectors in the approximation.

\begin{figure} \label{fig1}
\centering
\includegraphics[height=2.5in,width=3in]{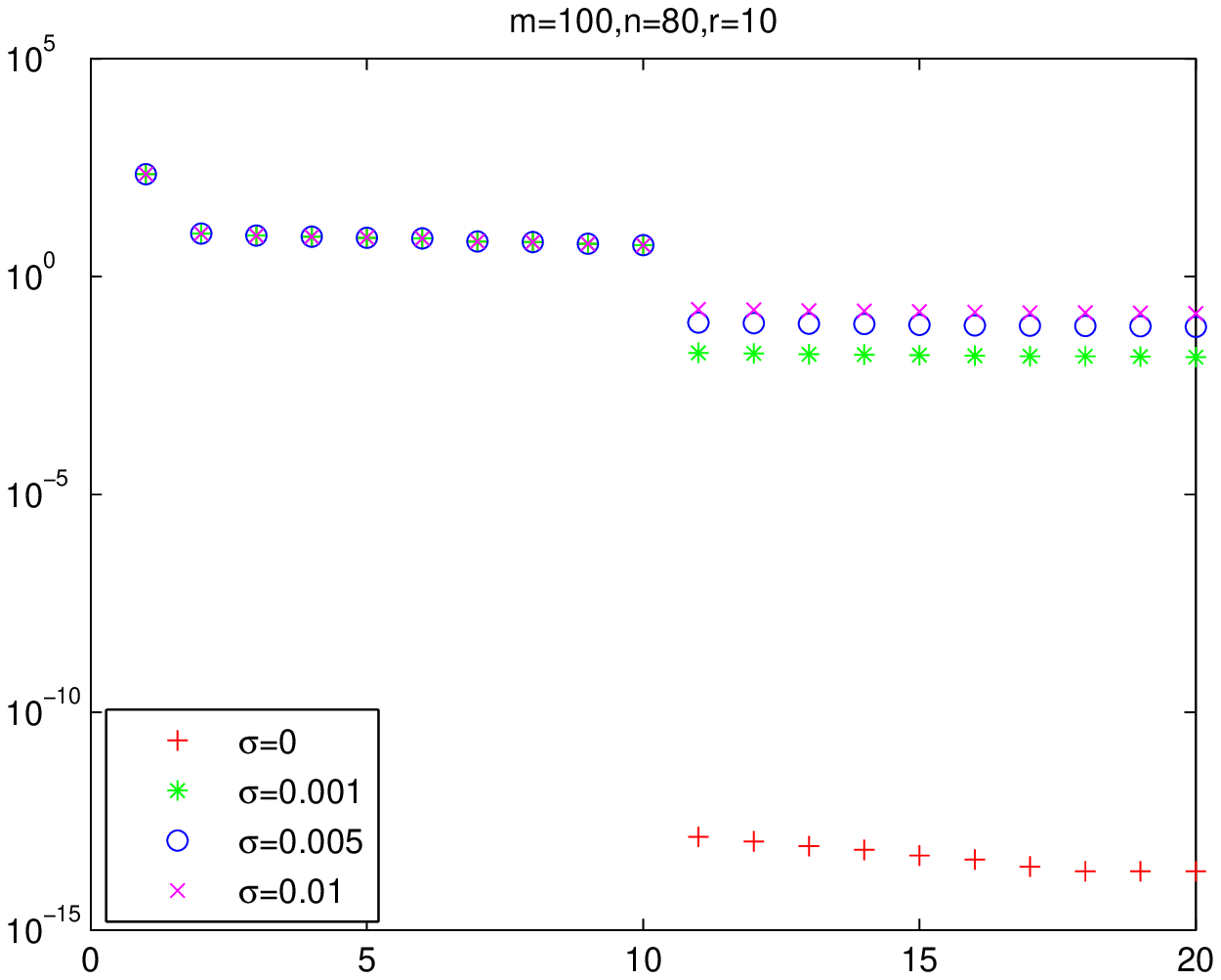}\includegraphics[height=2.5in,width=3in]{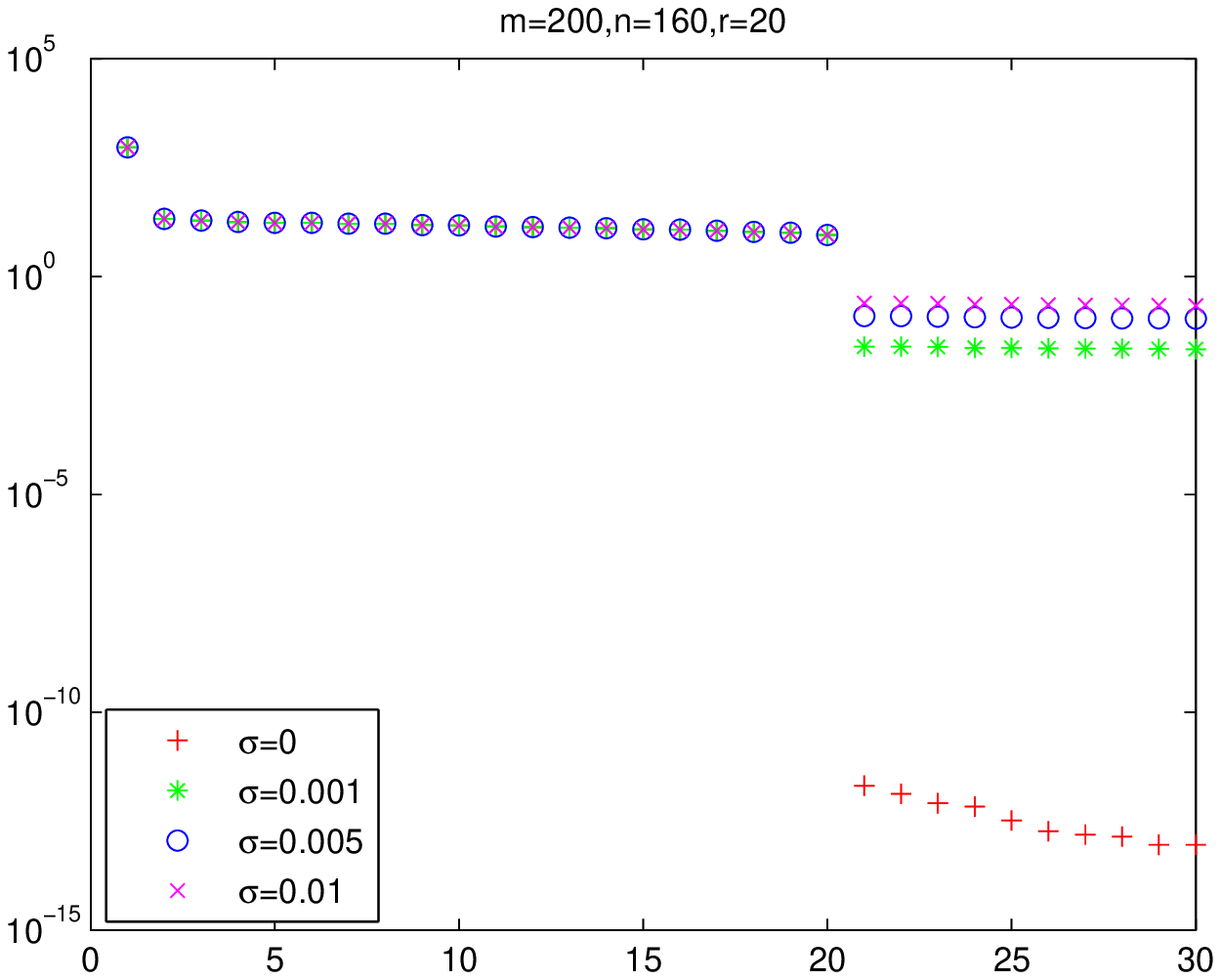}

\centerline{\quad \quad \quad \quad (a) \quad \quad \quad \quad \quad \quad \quad \quad \quad \quad \quad \quad \quad \quad \quad \quad\quad \quad \quad (b)}

\includegraphics[height=2.5in,width=3in]{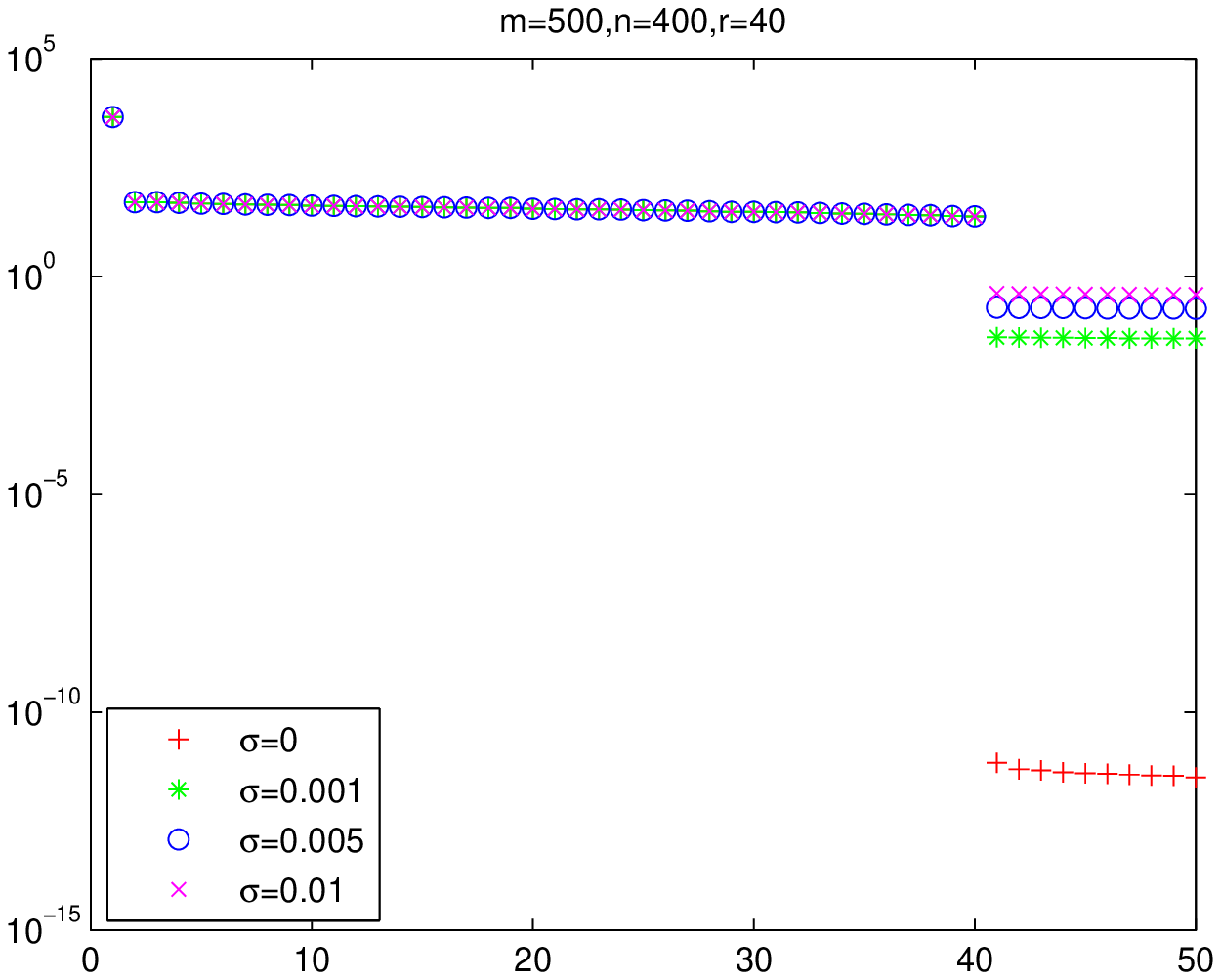}

\centerline{(c)}

\caption{The distribution of singular values of (a) the generated 100-by-80 matrix with actual rank 10,
(b) the generated 200-by-160 matrix actual rank 20, and (c) the generated 500-by-400 matrix with actual rank 40.}
\end{figure}

Similar to the first experiment, we randomly generated $m$-by-$n$ nonnegative matrices with actual ranks ($10,20,40$) and added
Gaussian noise of zero mean and variance $\sigma$.
In Figure 1, we display the distribution of singular values of the matrix approximation with rank $r$  ($=20,30,50$)
for the generated 100-by-80 matrix with actual rank 10, 200-by-160 matrix with actual rank 20, and
500-by-400 matrix with actual rank 40 respectively.
In this setting, we employ a matrix approximation with rank $r$ ($=20,30,50$) being larger than the actual rank ($=10,20,40$)
of the generated nonnegative matrix.
When there is no noise ($\sigma = 0$),
we see from Figures 1(a), 1(b), and 1(c) that there is a big jump in between the $k$-th singular value and the $(k+1)$th singular value
for actual ranks $k=10$, $k=20$, and $k=40$ respectively.
When there is a noise ($\sigma=0.001,0.005,0.01$), there is still a jump in between
the $k$-th singular value and the $(k+1)$th singular value ($k$ is the actual rank), and the height of the jump depends on the noise level.
This observation is also valid for other randomly generated matrices in the first experiment.
According to these figures, the distribution of singular values can provide us information to determine a suitable low rank
matrix approximation.

\begin{figure} \label{fig2}
\centering
\includegraphics[height=2.5in,width=3in]{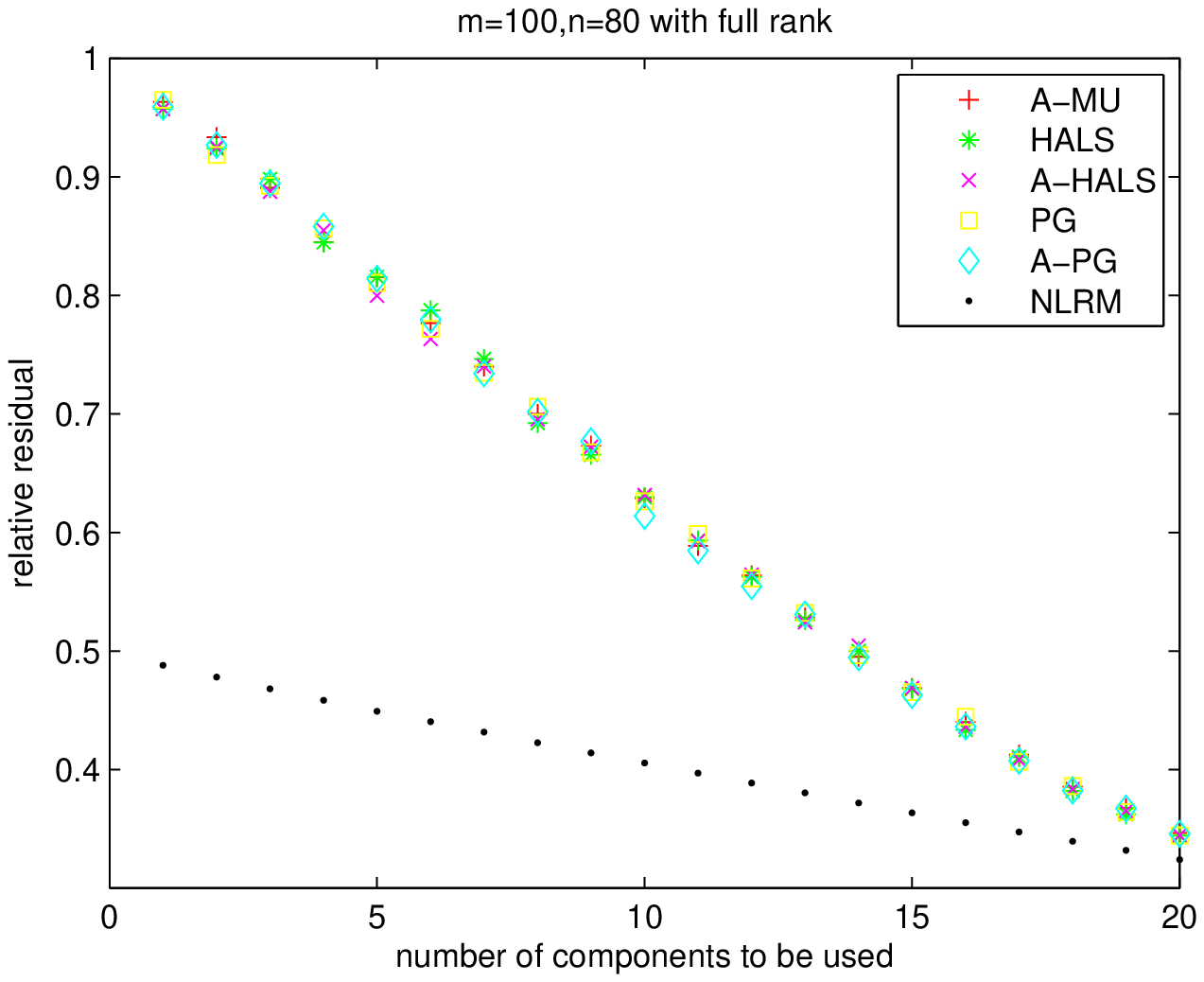}\includegraphics[height=2.5in,width=3in]{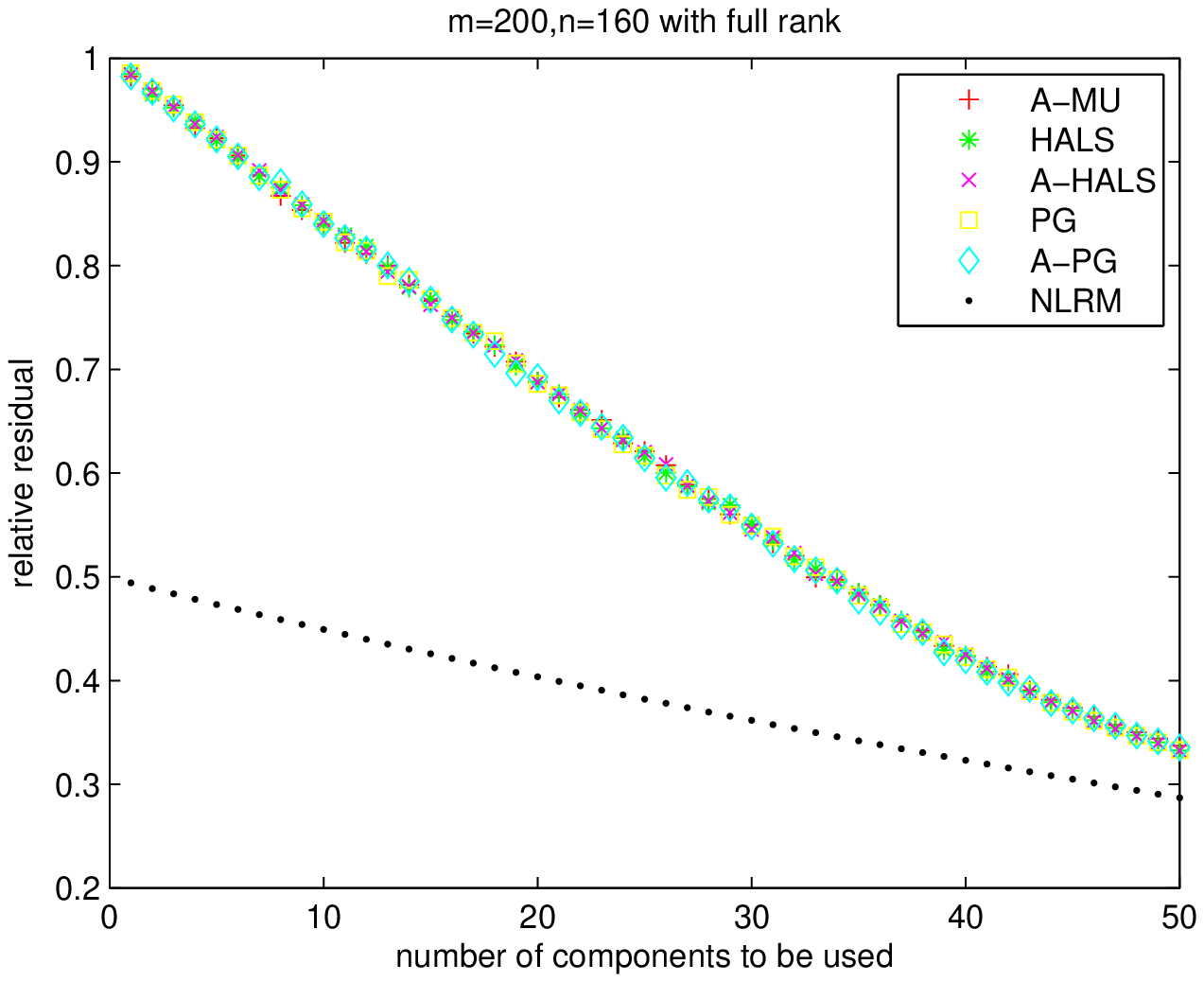}

\centerline{\quad \quad \quad \quad (a) \quad \quad \quad \quad \quad \quad \quad \quad \quad \quad \quad \quad \quad \quad \quad \quad\quad \quad \quad (b)}

\includegraphics[height=2.5in,width=3in]{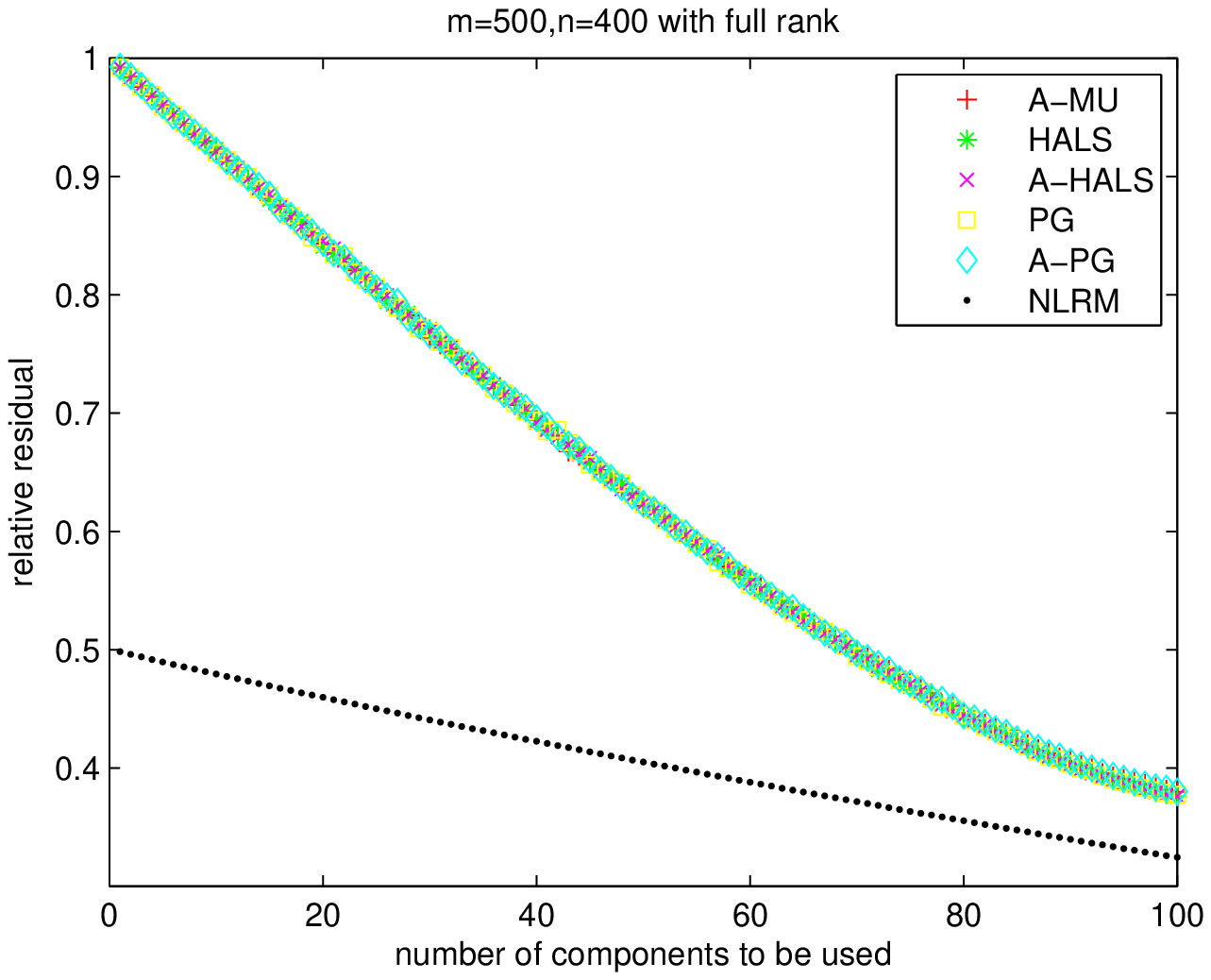}

\centerline{(c)}

\caption{The comparison of relative residuals with respect to the number of components to be used in the
matrix approximation with
(a) $r=20$ for the generated 100-by-80 matrix,
(b) $r=50$ for the generated 200-by-160 matrix, and
(c) $r=100$ for the generated 500-by-400 matrix.}
\end{figure}

\begin{figure} \label{fig3}
\centering
\includegraphics[height=2.5in,width=3in]{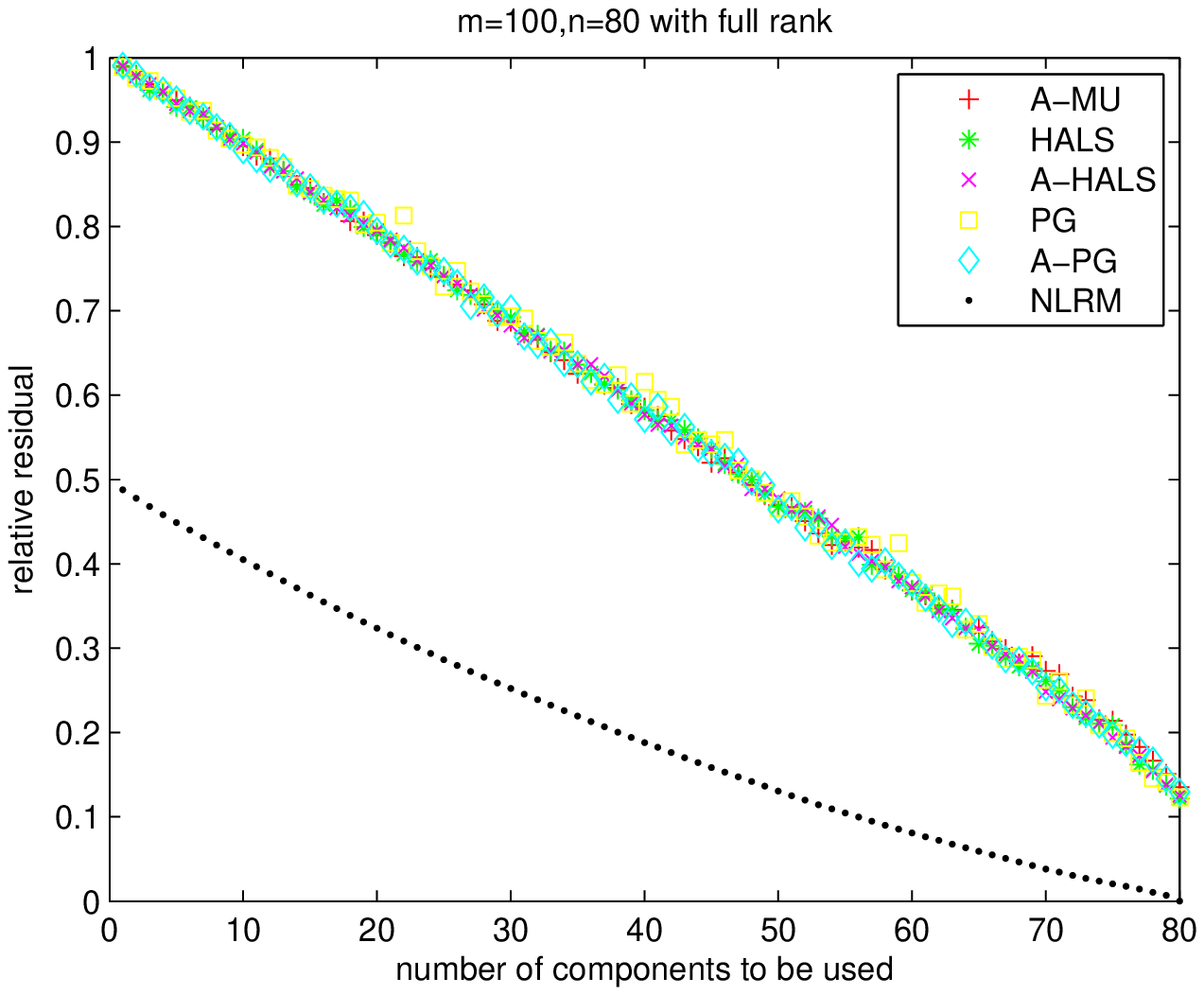}\includegraphics[height=2.5in,width=3in]{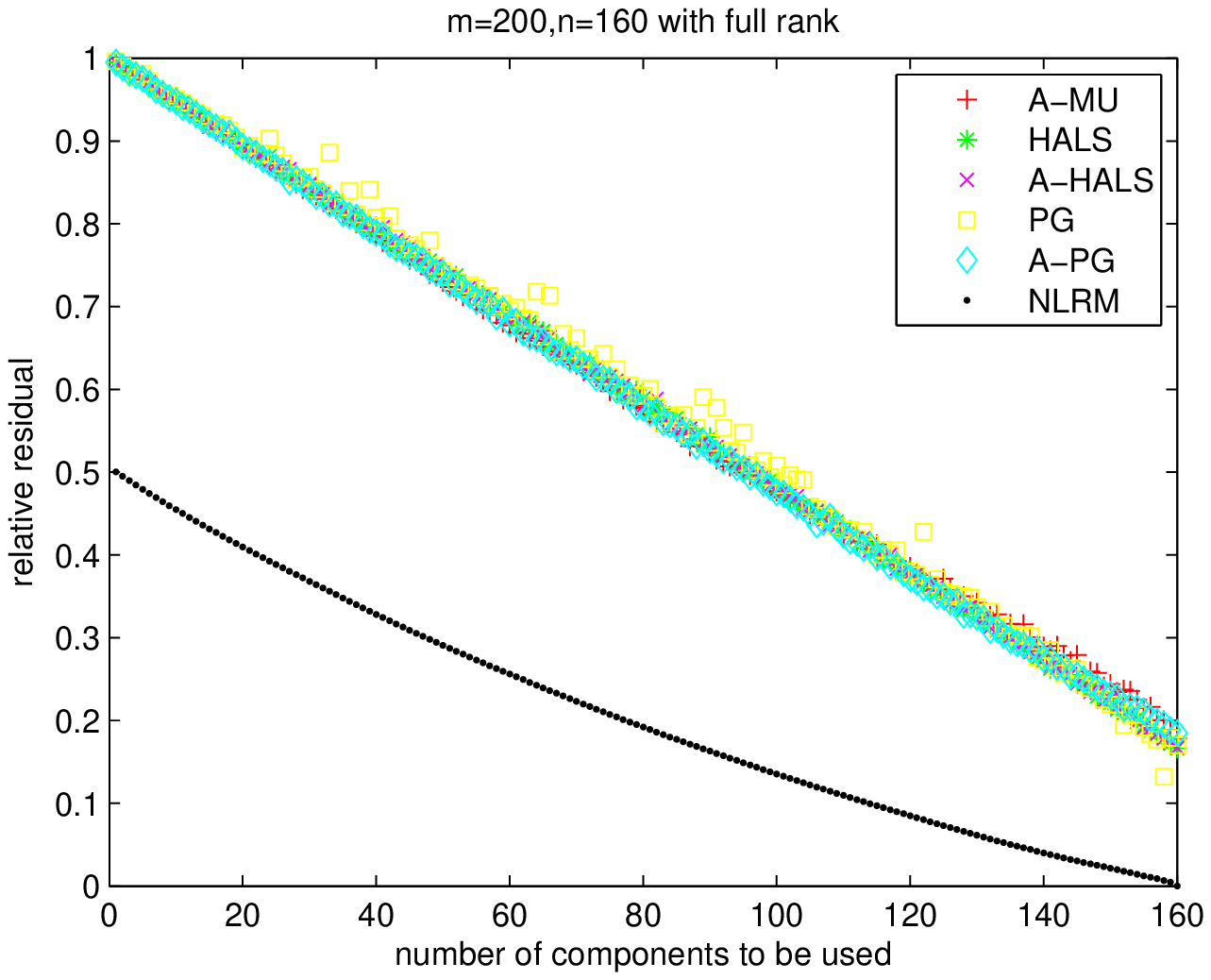}

\centerline{\quad \quad \quad \quad (a) \quad \quad \quad \quad \quad \quad \quad \quad \quad \quad \quad \quad \quad \quad \quad \quad\quad \quad \quad (b)}

\includegraphics[height=2.5in,width=3in]{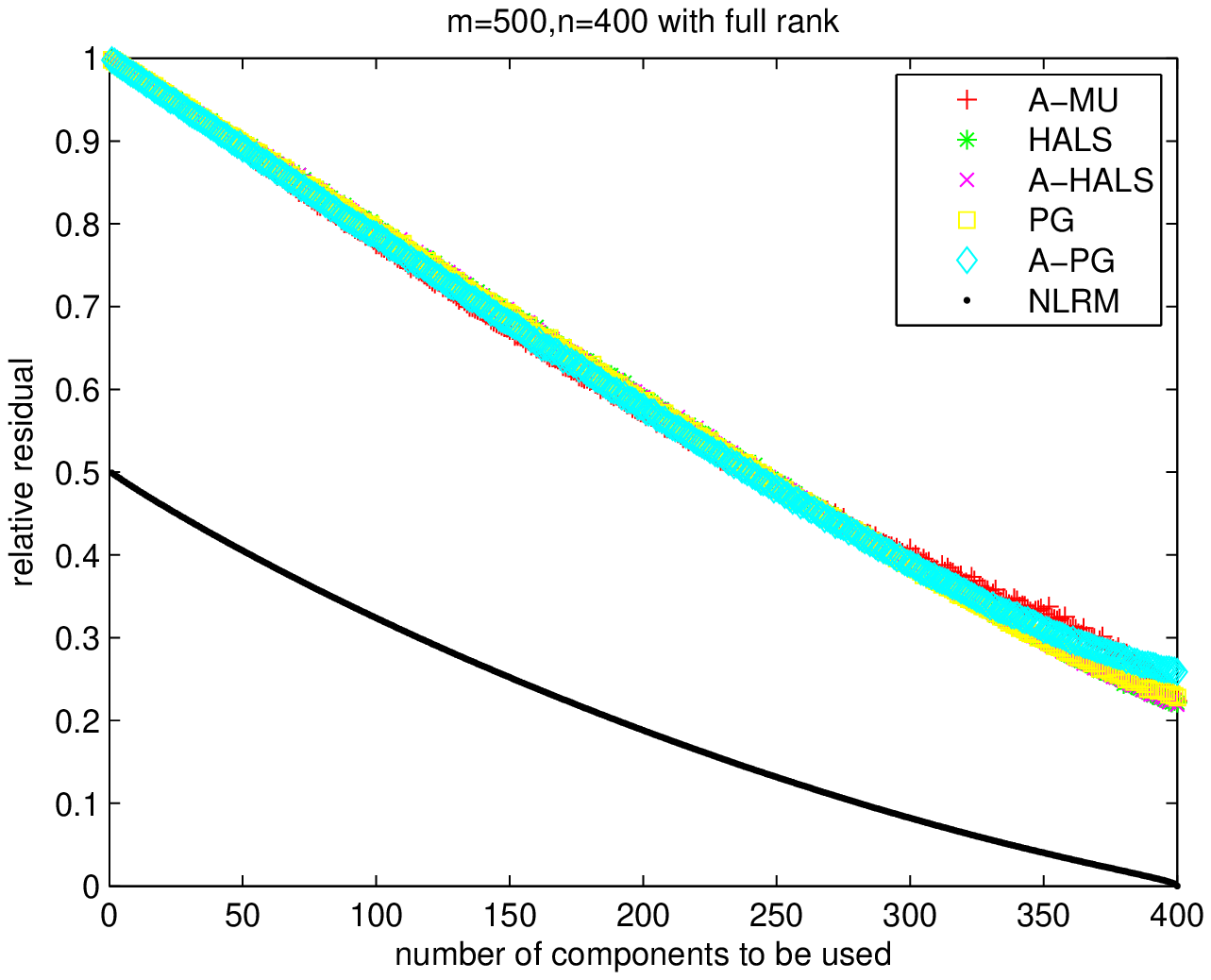}

\centerline{(c)}

\caption{The comparison of relative residuals with respect to the number of components to be used in the
matrix approximation with
(a) $r=80$ for the generated 100-by-80 matrix,
(b) $r=160$ for the generated 200-by-160 matrix, and
(c) $r=400$ for the generated 500-by-400 matrix.}
\end{figure}

On the other hand, we randomly generated $m$-by-$n$ nonnegative matrices with full rank.
In Figures 2-3, we display the relative residuals of the use of different numbers of the computed singular vectors by the proposed
NLRM algorithm. More precisely, the computed solution is given as follows:
$X_c = \sum_{i=1}^r \sigma_i u_i v_i^T$.
We plot $\| A - X_c(j) \|_F / \| A \|_F$ with respect to $j$ (the number of singular vectors to be used in the matrix approximation),
where $X_c(j) = \sum_{i=1}^{j} \sigma_i u_i v_i^T$
in the figures. Note that we employ singular vectors $u_i$ and $v_i^T$ according to the descending order
of the singular values $\sigma_i$ ($\sigma_1 \ge \sigma_2 \ge \cdots \ge \sigma_j$).
We see from the figures that when $j$ increases, the corresponding relative residual decreases.
These results show that the matrix approximation to the given nonnegative matrix
according to the ordering of singular values, is an effective strategy.
As a reference, when $r$ is equal to the full rank number, the relative residual
by the proposed NLRM algorithm, is close to the machine precision (around 1e-16).
In contrast, there is no index to show the columns of $m$-by-$r$ computed matrix $B_{comp}$ and the rows of $r$-by-$n$ computed
matrix $C_{comp}$ in the NMF approximation $B_{comp} C_{comp}$ to the nonnegative
matrix $A$. Here we normalize the row vectors of $C_{comp}$ such that the sum of squares of each row of $C_{comp}$ is equal to 1,
and then reorder the resulting column vectors of $B_{comp}$ according to their sum of squares.
We plot $\| A - X_{nmf}(j) \|_F / \| A \|_F$ with respect to $j$
where $X_{nmf}(j) = \sum_{i=1}^{j} [B_{comp}](:,i) [C_{comp}](i,:)$ and $[B_{comp}](:,i)$ is the $i$-th column of reordered $B_{comp}$ and
$[C_{comp}](i,:)$ is the $i$-th row of the normalized $C_{comp}$.
In Figures 2-3. In the figures, we see that when $j$ increases, the relative residual decreases by the testing NMF algorithms.
However, their relative residuals are still larger than those by the proposed NLRM algorithm.
It is interesting to note that
even $r$ is equal to the full rank number of the given nonnegative matrix, the relative residuals computed by the testing NMF algorithms are
not close to the machine precision. The reason may be that there is no guarantee that the testing NMF algorithms
give global optimal solutions.

Also we computed the relative residuals for face data matrix by the above mentioned procedure in Figure 4.
We see from the figure that the relative residuals by the proposed NLRM algorithm are significantly smaller than
those by the testing NMF algorithms for different values of $r=20,40,60,80,361$.
According to Figures 2-4, we summarize that
the proposed NLRM algorithm can provide
a significant index based on singular values that can be used to
identify important singular basis vectors in the approximation.

\begin{figure} \label{fig4}
\centering
\includegraphics[height=2.5in,width=3in]{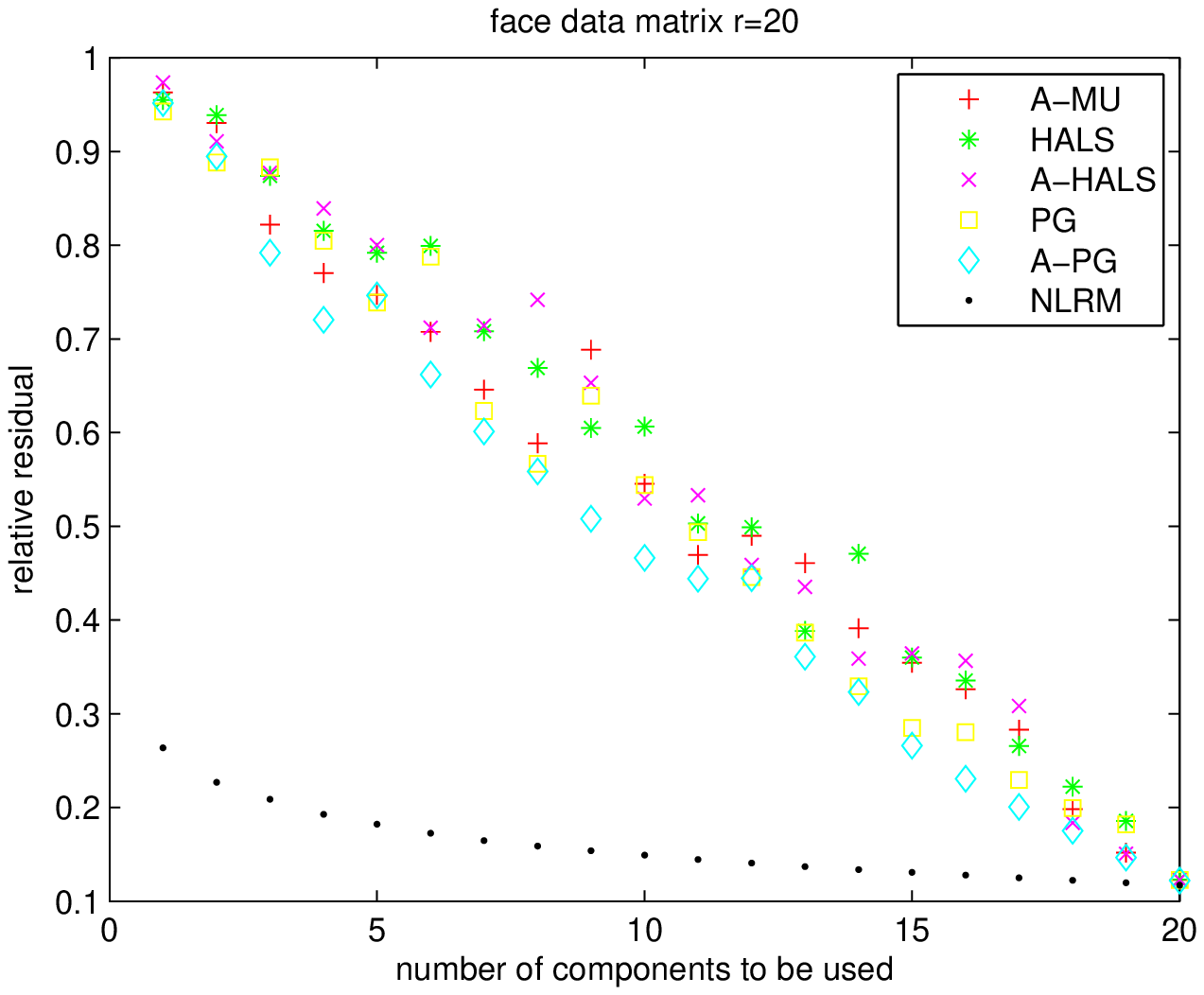}\includegraphics[height=2.5in,width=3in]{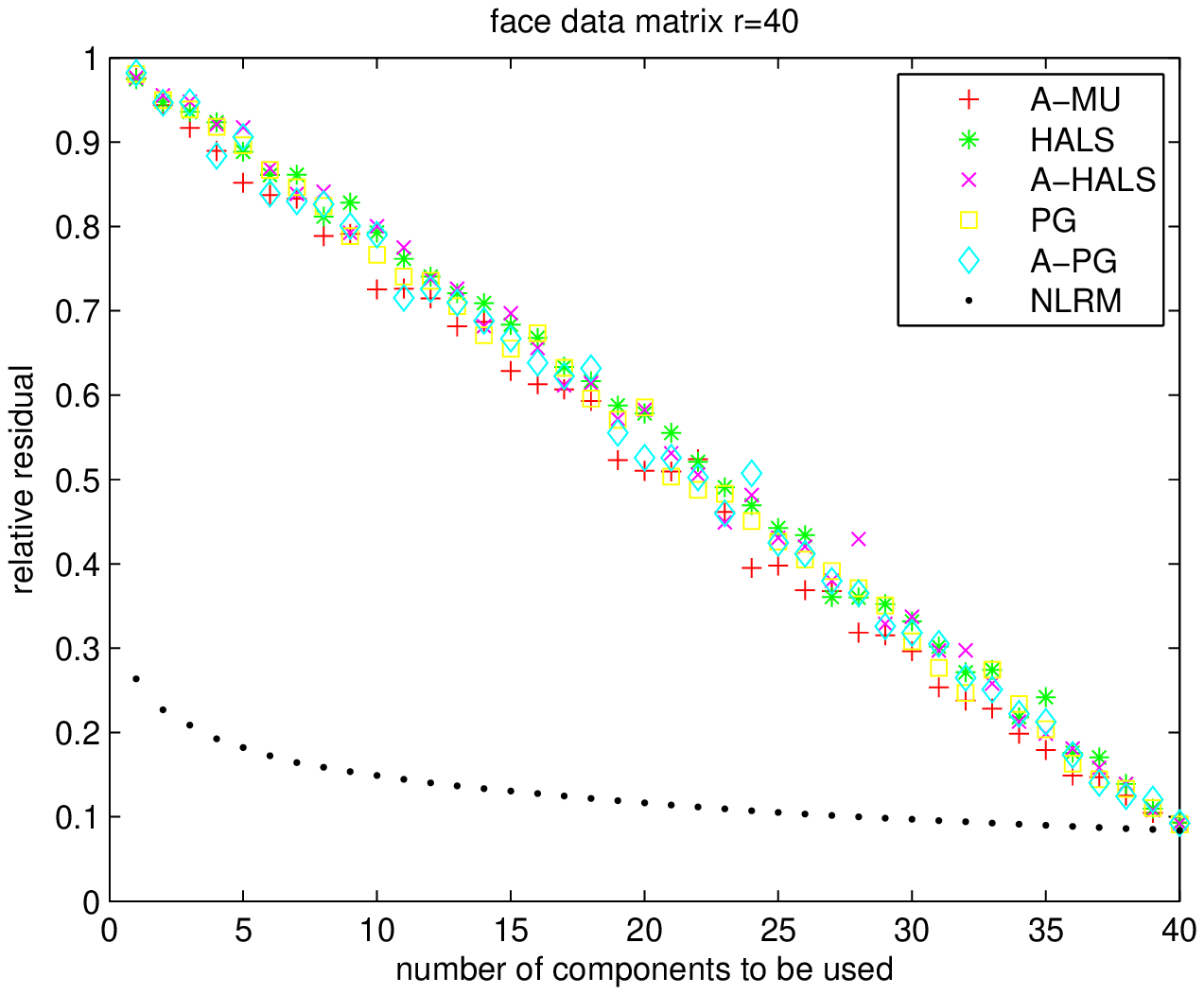}

\centerline{\quad \quad \quad \quad (a) \quad \quad \quad \quad \quad \quad \quad \quad \quad \quad \quad \quad \quad \quad \quad \quad\quad \quad \quad (b)}

\includegraphics[height=2.5in,width=3in]{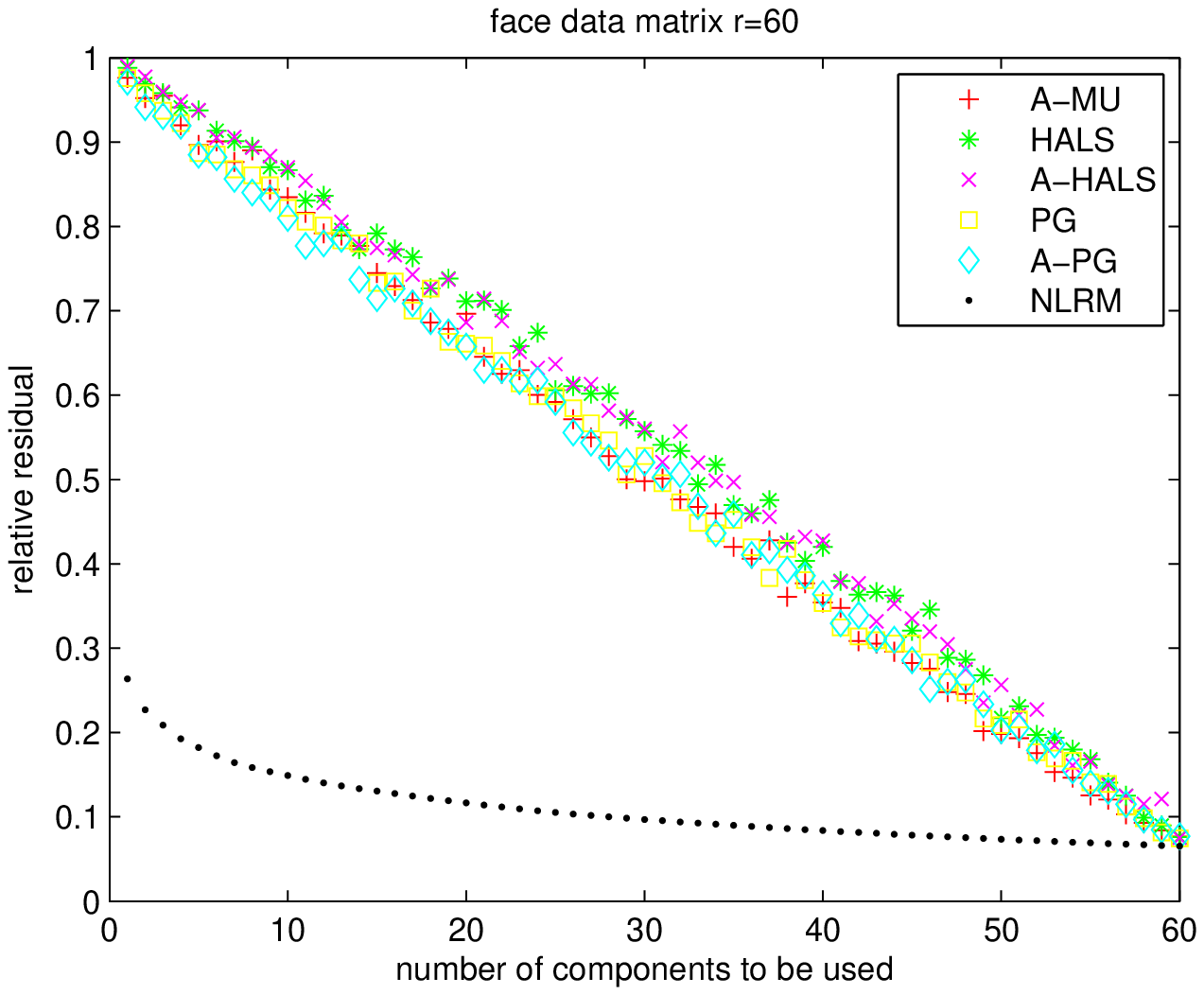}\includegraphics[height=2.5in,width=3in]{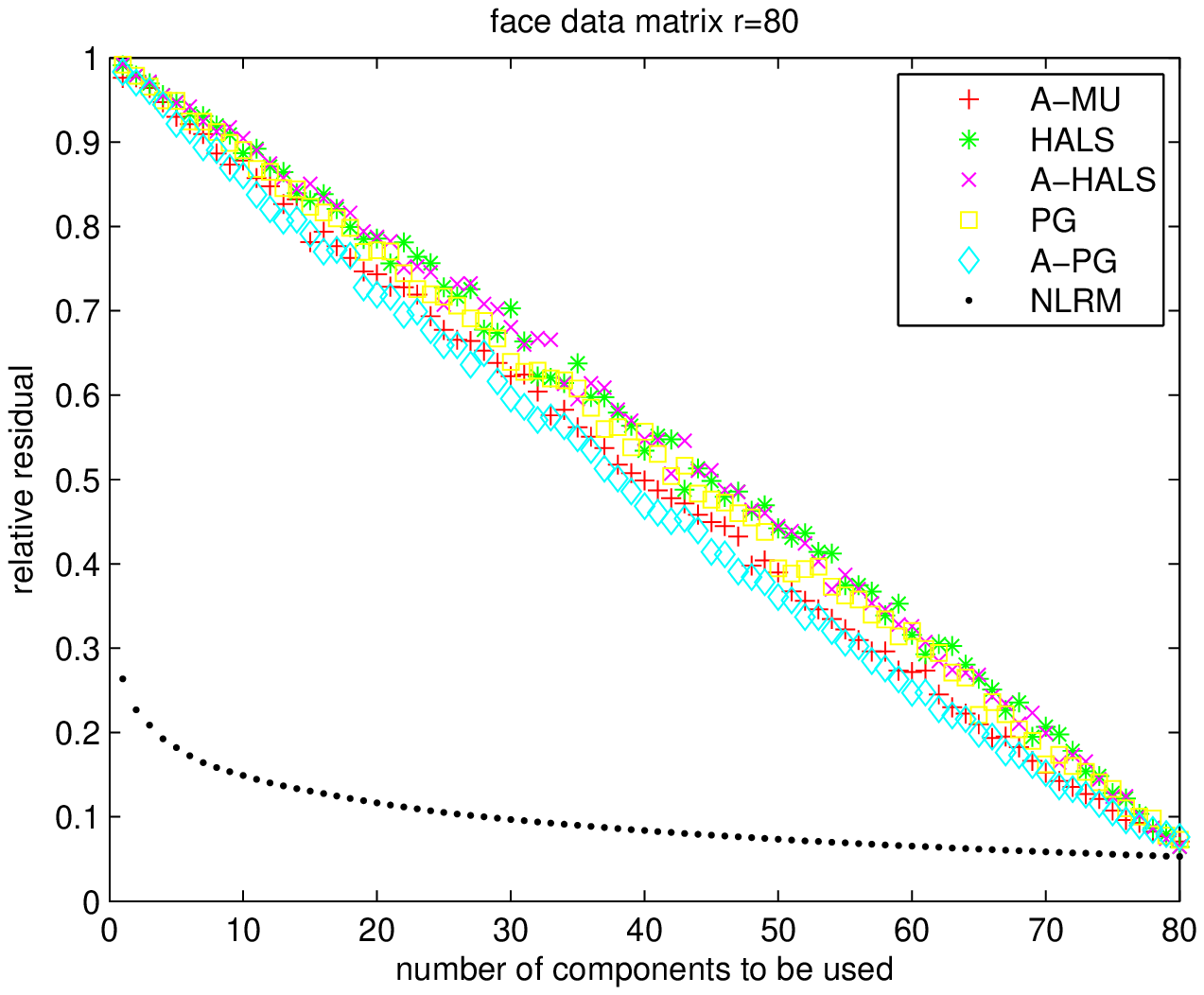}

\centerline{\quad \quad \quad \quad (c) \quad \quad \quad \quad \quad \quad \quad \quad \quad \quad \quad \quad \quad \quad \quad \quad\quad \quad \quad (d)}

\includegraphics[height=2.5in,width=3in]{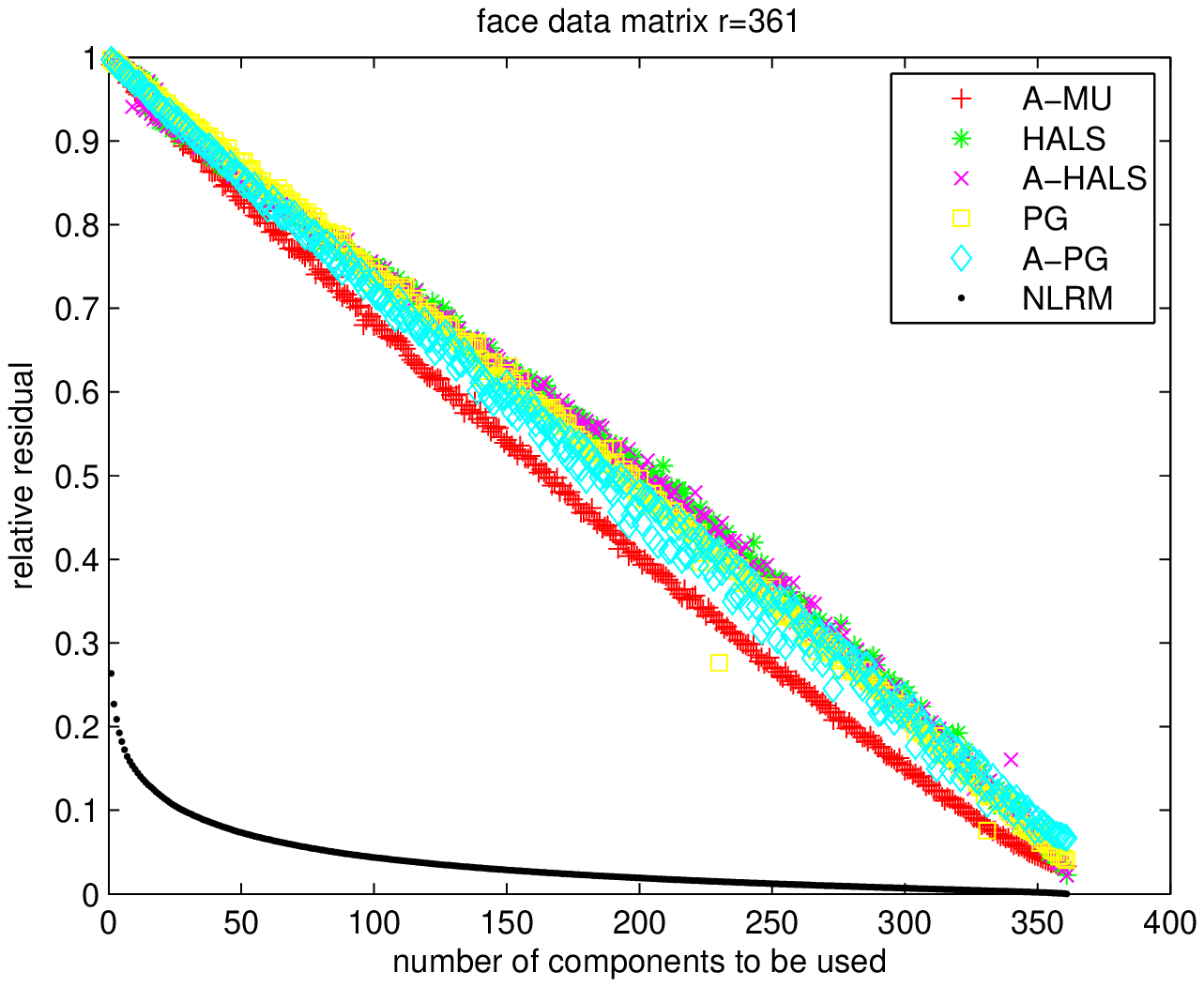}

\centerline{(e)}

\caption{The comparison of relative residuals with respect to the number of components to be used in the
face data matrix approximation with (a) $r=20$ (b) $r=40$ (c) $r=60$ (d) $r=80$ (e) $r=361$ (full rank).}
\end{figure}

\section{Concluding Remarks}

In this paper, we proposed and developed a new algorithm for computing NLRM approximation for
nonnegative matrices. The new method is different from classical nonnegative matrix
factorization method. We have shown the convergence of the proposed algorithm based on the results in manifold.
Moreover, we have demonstrated numerical results that
the minimized distance by the proposed NLRM method can be smaller than
that by the NMF method. According to the ordering of singular values,
the proposed method identifies important singular basis vectors, while this information cannot be
obtained in the classical NMF. As a future research work, there are several areas to be
studied:
\begin{enumerate}
\item
The computational cost of the proposed
method involves the computation of singular value decomposition. We plan to study tangent space
method for low rank matrix projection so that the computational complexity can be reduced for
large scale data science applications.
\item
We study applications of
nonnegative low rank matrix approximation and check how the proposed method can
provide a significant index based on singular values that can be used to
identify important singular basis vectors in the approximation.
\item
In classical NMF applications, researchers have suggested to use the other norms (such as $l_1$ norm and KL divergence)
in data fitting
$\| A - B C \|$ instead of Frobenius norm to deal with other machine learning applications.
It is interesting to develop the related algorithms for nonnegative low rank matrix approximation.
\end{enumerate}

%\bibliography{nonnegative}

\section*{Supplementary Information}

The supplementary information consists of the theoretical proofs
of Theorems \ref{th1}-\ref{thm_convergence}.

Some basic definitions, propositions of algebraic geometry and differential geometry are needed to prove the main results. We only provide some results here and for
details we refer to \cite{absil2009optimization,lee2013smooth} and references therein.

In order to show the convergence of  Algorithm \ref{ag1}, we need to prove the intersection of
$\mathcal{M}_{r}$ given as  \eqref{v1} and $\mathcal{M}_{n}$ given as \eqref{v2} is a manifold first.
In many real applications, some manifolds are actually real algebraic varieties which are defined as the vanishing of a set of polynomials on $\mathbb{R}^{n}$, thus some algebraic  geometry methods can help us to study the above problem. Note that all complex varieties are real varieties, but not conversely. For a given real algebraic variety $\mathcal{V}\in \mathbb{R}^{n}$, if we identity $\mathbb{R}^{n}$ as a subset of $\mathbb{C}^{n}$ and denote $ \mathbb{I}_{\mathbb{R}}(\mathcal{V}) $ as the set of real polynomials that vanish on $\mathcal{V}$, then $\mathcal{V}$ has a related complex variety given by its Zariski closure
$$
\mathcal{V}_{Zar} = \{ z \in \mathbb{C}^{n} : p(z)=0, \ \forall p \in \mathbb{I}_{\mathbb{R}}(\mathcal{V}) \},
$$
which is defined as the subset in $\mathbb{C}^{n}$ of common zeros to all polynomials that vanish on $\mathcal{V}$. Moreover,
 H. Whitney in \cite{whiteneyvarieties} showed that a real algebraic variety $\mathcal{V}$ can be decomposed as $\mathcal{V}=\bigcup_{j=0}^{m} \mathcal{V'}_{j}$ where each $\mathcal{V'}_{j}$ is either void or a $\mathbb{C}^{(\infty)}$-manifold with dimension $j$. If $\mathcal{V'}_{m}\neq \emptyset$,  then $m$ equals the algebraic dimension of $\mathcal{V}_{Zar}$. Each $\mathcal{V'}_{j}$ contains at most a finite number of connected components.  This result shows that the main part of a variety is a manifold.

Varieties have singular points which is different from manifolds. Note that a point $A\in \mathcal{V}$ is nonsingular if it is nonsingular in the sense of algebraic geometry as an element of $\mathcal{V}_{zar}$. This result changes to be much simpler, if we restrict on the irreducible variety.  Here a algebraic variety $\mathcal{V}$ is said to be irreducible if there does not exist non-trivial real algebraic varieties $\mathcal{V}'_{1}$
and $\mathcal{V}'_{2}$, such that $\mathcal{V}=\mathcal{V}'_{1}\bigcup \mathcal{V}'_{2}$.
 Denote $\nabla$ as the gradient operator and set
${\cal N}_{\cal V}(z) = \{ \nabla p(z) : p \in \mathbb{I}_{\mathbb{R}}(\mathcal{V}) \}$. Suppose that $\mathcal{V}\in \mathbb{R}^{n}$ is an irreducible real algebraic variety of dimension $m$, then $\dim {\cal N}_{\mathcal{V}}(z)\leq n-m$ and $z\in \mathcal{V}$ is non-singular if and only if $\dim {\cal N}_{\mathcal{V}}(z)=n-m$.  In  practice, it is not easy to check a given variety is irreducible or not,  thus the following results are needed.

\begin{definition}[Definition 6.7 in \cite{andersson2013alternating}]\label{newdf}
  Suppose we are given a number $j\in \mathbb{N}$ and an index set $I$ such that for each $i\in I$, there exist an open connected $\Omega_{i}\subseteq \mathbb{R}^{j}$ and a real analytic map $\phi_{i}:\Omega_{i}\rightarrow \mathcal{V}.$  Then $\mathcal{V}$ is said to be covered with analytic patches, if for each $A\in \mathcal{V},$ there exists an $i\in I$ and a radius $r_{A}$ such that
$$\mathcal{V}_{rn}\cap Ball_{\mathbb{R}^{n}}(A,r_{A})=Im \phi_{i}\cap  Ball_{\mathbb{R}^{n}}(A,r_{A}).$$
\end{definition}

\begin{proposition}[Proposition 6.8 in \cite{andersson2013alternating}]\label{pr1}
Let V be a real algebraic variety. If $\mathcal{V}$ is connected and can be covered with analytic patches, then $\mathcal{V}$ is irreducible.
 %The same conclusion holds if a dense subset of $\mathcal{V}$ can be given as the image of one real analytic function $\phi.$
\end{proposition}

Moreover, the following proposition proposed us a simple method to  compute the dimensions of some varieties involved which are tricky to compute.

\begin{proposition}[Proposition 6.9 in \cite{andersson2013alternating}]\label{pr2}
Under the assumption of Proposition \ref{pr1}, suppose in addition that an open subset of $\mathcal{V}$ is the image of a bijective real analytic map defined on a subset of $\mathbb{R}^{d}$. Then $\mathcal{V}$ has dimension $d$.
\end{proposition}

With the above tools in hand we can prove the following results. The idea of the proof follows from [2].

\vspace{3mm}
\noindent
{\bf Proof of Theorem \ref{th1}}:
The proof can be divided into two parts. Firstly, we need to prove the following set
\begin{align}\label{v5}
\mathcal{V}_{rn}=\left\{X\in \mathbb{R}^{m\times n}, ~ \rank(X)\leq r,~X_{i,j}\geq 0,~ i=1,...,m,j=1,...,n\right\}
\end{align}
is an irreducible variety with dimension $(m+n)r-r^2$. Then we need to prove $\mathcal{V}_{rn}^{ns}$
(the set of all the non-singular points in $\mathcal{V}_{rn}$) equals the set of all the nonnegative matrices with rank equal to $r$.

%In the first step,
%we mainly prove $\mathcal{V}_{rn}$ is an irreducible variety with dimension $(m+n)r-r^2$.
Denote $\mathcal{K}$ as set of $m\times n$ matrices over $\mathbb{R}$. If all elements of a matrix $A\in\mathcal{K}$ are considered as variables, $\mathcal{K}$ is a linear manifold with dimension $mn$.
The same statement is also satisfied for the nonnegative matrix set $\mathcal{M}_{n}$.
Denote
$\mathcal{V}_{r}:=\left\{ X\in \mathbb{R}^{m\times n}, \ \rank(X)\leq r\right\}$.
Recall that a matrix in $\mathcal{K}$ has rank being greater than $r$ if and only if one can find a nonzero $(r+1)\times (r+1)$ invertible minor. The determinant of each such minor is a polynomial, and $\mathcal{V}_{r}$ is clear the variety obtained from the collection of such polynomials. Thus $\mathcal{V}_{r}$ is a real algebraic variety. The same statement is also true for $\mathcal{V}_{rn}$, which is
derived by adding that a matrix entry is the square of a variable (i.e., nonnegative constraints)
to the variety defining on $\mathcal{V}_{r}$.
In order to use Proposition \ref{pr1}, we need to show that $\mathcal{V}_{rn}$ is covered with analytic patches  and connected,
and then $\mathcal{V}_{rn}$ can be irreducible.

\noindent
(i) (Covered with analytic patches)
%Let $A\in \mathcal{V}_{rn}$.
%It follows from the singular value decomposition theory \cite{golub2012matrix} that there exist two column orthogonal matrices $U\in %\mathbb{R}^{m\times r'}$, $V\in \mathbb{R}^{n\times r'}$ and a diagonal matrix $\Sigma\in \mathbb{R}^{r'\times r'}$
%(containing the singular values of $A$) such that
%\begin{align}\label{xin1}
%A=U \Sigma  V^{T}.
%\end{align}
%Here $r'$ is equal to the rank of $A$.
%Note that some main diagonal entries of $\Sigma$ can be zero as $A$ has rank being less than or equal to $r$.
%A matrix written as in \eqref{xin1} has rank less or equal to $r$,
%therefore $\mathcal{V}_{rn}$ can be covered by a real polynomial.
%In order to use Proposition \ref{pr1}, we need to show that $\mathcal{V}_{rn}$ is connected and can be covered with analytic patches,
%and then $\mathcal{V}_{rn}$ can be irreducible.
%\noindent
%(i) (Covered with analytic patches)
%Let $A\in \mathcal{V}_{rn}$.
%We express $A$ as in (\ref{xin1}).
%Suppose that $\pi(i)=r'$ for all $i=1,...,m$.
Denote $U\in \mathbb{R}^{m\times r'}$, $V\in \mathbb{R}^{n\times r'}$ and a diagonal matrix $\Sigma\in \mathbb{R}^{r'\times r'}$ where
$r' \le r$. Here we assume $r' \ge 1$. For $r' =0$, it is a trivial case.
Suppose that $\pi(i)=r'$ for all $i=1,...,m$.
We set $u_{i,\pi(i)}$ as an undetermined variable.
Now we can construct a real analytic mapping $\theta$ from $\hat{U}$ (the first $(r'-1)$ rows of
$U$),  $\Sigma$ and  $V$  to $\mathcal{V}_{rn}$ as follows:
\begin{eqnarray} \label{ineq1}
& &
\theta(u_{1,1},u_{1,2},...,u_{1,r'-1},...,u_{m,1},...,u_{m,r'-1},\sigma_{1},...,\sigma_{r'},v_{1,1},...,v_{n,r'}) \nonumber \\
& =&
 \left(
             \begin{array}{cccc}
               u_{1,1} & \cdots & u_{1,r'-1} & u_{1,r'} \\
               u_{2,1} & \cdots & u_{2,r'-1} & u_{2,r'} \\
              \vdots & \ddots& \vdots & \vdots \\
               u_{m,1} & \cdots & u_{m,r'-1} & u_{m,r'}\\
             \end{array}
           \right)\left(
                    \begin{array}{cccc}
                      \sigma_{1} & 0 & \cdots & 0 \\
                      0 & \sigma_{2} & \cdots & 0 \\
                      \vdots & \vdots &\ddots & \vdots \\
                      0 & 0 & \cdots & \sigma_{r'} \\
                    \end{array}
                  \right)\left(
             \begin{array}{cccc}
               v_{1,1} & \cdots & v_{1,r'-1} & v_{1,r'} \\
               v_{2,1} & \cdots & v_{2,r'-1} & v_{2,r'} \\
              \vdots & \ddots& \vdots & \vdots \\
               v_{n,1} & \cdots & v_{n,r'-1} & v_{n,r'}\\
             \end{array}
           \right)^{T}. \nonumber \\
\end{eqnarray}
Note that the entries of the matrix in (\ref{ineq1}) can be nonnegative when
the associated inequalities of $m$ undetermined variables are satisfied.
Such $m$ undetermined variables are decoupled in these inequalities.
There are infinitely many real solutions of $m$ undetermined variables for given values of $\hat{U}$, $\Sigma$, and $V$.
It is saying that $\mathcal{V}_{rn}$ is the image of $\theta$.
Let $\Gamma$ be a particular connected component of $(\hat{U}, \Sigma, V)$.
We establish a function $\psi$ with $\pi$ and $\Gamma$ as follows:
\begin{equation}\label{tr1}
\psi_{\pi,\Gamma}(y) = U(y) \Sigma(y) V(y)^T, ~~y \in \Gamma.
\end{equation}
Denote $\mathbb{I}$ as the set of all possible $\pi$ and $\Gamma$.
It can be found that for each matrix in $\mathcal{V}_{rn}$ is in the image of at least one $\psi_{\pi,\Gamma}$ where
$(\pi, \Gamma) \in \mathbb{I}$.
Then by Definition \ref{newdf},
$\mathcal{V}_{rn}$ can be covered by $\{ \psi_{\pi,\Gamma} \}_{(\pi, \Gamma) \in \mathbb{I}}$.

\noindent
(ii) (Connected)
In order to show $\mathcal{V}_{rn}$ is connected, we need to show that for any two nonnegative matrices $A, B\in \mathcal{V}_{rn}$,
there exist a continuous map $f$ from the unit interval $[0,1]$ to $\mathcal{V}_{rn}$ such that $f(0) = A$ and $f(1) = B$.
Without loss of generality, we show that
for an arbitrary $A \in \mathcal{V}_{rn}$, it is connected with $\textbf{1}$ matrix (all the entries are equal to 1) instead.
It is sufficient to prove $\mathcal{V}_{rn}$ is path connected.
Suppose that $A, B\in \mathcal{V}_{rn}$ are arbitrary and path connected with the $\textbf{1}$ matrix, respectively. Thus there are two continuous  maps $f$ and $g$  which are from the unit interval $[0,1]$ to $\mathcal{V}_{rn}$  with $f(0) = A$, $f(1) = \textbf{1}$,
$g(0)=\textbf{1}$ and $g(1)=B$. Setting $\tau(x)=(1-x)f(x)+xg(x)$, it is easy to see that $\tau(x)$ is continuous and satisfying $\tau(0)=f(0)=A$ and $\tau(1)=g(1)=B$. Then $\mathcal{V}_{rn}$ is  path connected.
Let $A \in\mathcal{V}_{rn}$. Suppose that all the singular values of $A$ are ordered decreasingly with
$\sigma_{1}=1$ (if $\sigma_{1} \neq 1$, we can divide or multiply some factor such that $\sigma_{1}=1$.
Set  $\pi(i)=r'=rank(A)$, $i=1,...,m$ as above and choose
$\Gamma$ such that the matrix representation can be expressed as \eqref{tr1}.
Now if $\sigma_{2} \neq 1$, we
can proceed until it is not without leaving $\Gamma$. Then the values of $y$ corresponding to the first and second column of
$\hat{U}$ and $V$ can be continuously moved until all elements of the first column of $\hat{U}$ and $V$ are positive, respectively.
At this point, we can
reduce all values of $\hat{U}$ and $V$ except the first column to zero, increase the first value of each row whenever necessary to stay in
$\Gamma$. Then we can move $y$ so that the values in the first column become the same. Finally, we can let these values increase simultaneously until
they reach $1$. Then the matrix $\mathbf{1}$ is derived, which is saying that $\mathcal{V}_{rn}$ is connected.
Hence it is irreducible.

\noindent
(iii) Now we would like to apply Proposition \ref{pr2} so that
$\mathcal{V}_{rn}$ possesses dimension $(m+n)r-r^2$.
Here we need to find an open subset of $\mathcal{V}_{rn}$ which can be expressed as the image of a bijective real analytic map
defined on a subset of $\mathbb{R}^{(m+n)r-r^2}$.
Now we set $\Upsilon$ as a subset of $\mathcal{V}_{rn}$ in which the matrix rank is $r$.
There are $r$ singular values.
According to (\ref{ineq1}) and the singular value decomposition theory \cite{golub2012matrix},
we know that for a matrix in $\Upsilon$, the columns of $\hat{U}$ are orthogonal and the dimension is $(mr-\frac{r(r+1)}{2}-m)$,
and the rows of $V$ are orthogonal and the dimension is $(nr-\frac{r(r+1)}{2})$.
In addition, there are $m$ undetermined variables in the last column of $U$,
and there are $r$ independent variables on the main diagonal matrix $\Sigma$.
Therefore, the combined dimension is $(m+n)r-r^2$.
Now we can identify three sets of matrices with
$\mathbb{R}^{mr-\frac{r(r+1)}{2}}$, $\mathbb{R}^{nr-\frac{r(r+1)}{2}}$ and $\mathbb{R}^{r}$ respectively.
Denote the inverses of these identifications by
\begin{align}\label{pj1}
\iota_{1} : \mathbb{R}^{mr-\frac{r(r+1)}{2}}\rightarrow \mathbb{R}^{m\times r};~ \iota_{2} : \mathbb{R}^{r}\rightarrow \mathbb{R}^{r\times r};~\iota_{3} : \mathbb{R}^{mr-\frac{r(r+1)}{2}}\rightarrow \mathbb{R}^{n\times r};
 \end{align}
and denote $\Omega \subset \mathbb{R}^{(m+n)r-r^2}$ as the open set corresponding to those matrices having the same structure as $\Upsilon$.
Define $\phi: \Omega \rightarrow \mathcal{V}_{rn}$ by
\begin{equation}\label{qj2}
\phi(y_{1},y_{2},y_{3})=\iota_{1}(y_{1})\times \iota_{2}(y_{2})\times \iota_{3}(y_{3})^{T}.
\end{equation}
It is easy to see that $\phi$ is bijective correspondence with an open set of $\mathcal{V}_{rn}$ and $\phi$ is a polynomial. By
Proposition \ref{pr2}, $\mathcal{V}_{rn}$ possesses dimension $(m+n)r-r^2$ as desired.

In the second step, we mainly prove $\mathcal{V}_{rn}^{ns}$ equals all the nonnegative matrices with rank equal to $r$.
Recall  that $\mathcal{V}_{rn}$  is an irreducible real algebraic variety with dimension $(m+n)r-r^2$, then we need to prove
\begin{align*}
\dim \mathcal{N}_{\mathcal{V}_{rn}}(A)=mn-(m+n)r+r^2,
\end{align*}
if and only if $\rank(A)=r$. It is sufficient  to show that $\dim \mathcal{N}_{\mathcal{V}_{rn}}(A)\leq mn-(m+n)r+r^2$ is strict
when $\rank(A)<r$ and  the reverse inequality holds when $rank(A)=r$.
%Define $\omega: \mathbb{R}^{mn}\rightarrow \mathcal{K}$.
%If $p$ is a polynomial on $\mathcal{K}$,
%then write $\nabla p$ in stead of the right form $\omega(\nabla(p\circ \omega))$.
Now for a given  polynomial $p\in \mathbb{I}(\mathcal{V}_{rn})$, and
two orthogonal matrices $U$ and $V$, $q_{U,V}(\diamond)=p(U\diamond V^{T})$,
which is also in $\mathbb{I}(\mathcal{V}_{rn})$. Let $A$ be  a fixed nonnegative matrix with $\text{rank}(A)=k\leq r$.
By using the singular value decomposition, we can produce two orthogonal matrices $U$ and $V$ such that $$UAV^{T}=\diag\{\sigma_{1},...,\sigma_{k},0,...,0\}=E_{k},$$ where the first $k$
diagonal values of $E_{k}$ are positive numbers and $0$ otherwise. In particular $\nabla q_{U,V}(A)=U\nabla p(UAV^{T})V^{T}=U\nabla p(E_{k})V^{T}$, implying that
$\dim \mathcal{N}_{\mathcal{V}_{rn}} (A)=\dim \mathcal{N}_{\mathcal{V}_{rn}} (E_{k})$. If $k=r,$   all the $(r+1)\times (r+1)$ subdeterminants of the  matrices in $\mathcal{K}$ form polynomials in $\mathbb{I}(\mathcal{V}_{rn})$, thus
\begin{align*}
\dim \mathcal{N}_{\mathcal{V}_{rn}} (E_{r})\geq mn-(m+n)r+r^2,
\end{align*}
which proves that any rank $r$ element of $\mathcal{V}_{rn}$ is nonsingular. Moreover, if $k< r$, consider some fixed $u\in \mathbb{R}_{+}^{m}$ and $v\in \mathbb{R}_{+}^{n}$, define
the map $\vartheta_{u,v}:\mathbb{R_+}\rightarrow \mathcal{V}_{rn}$ via $\vartheta_{u,v}(x)=E_{k}+xuv^{T}$.  Considering
\begin{align*}
\dim\left\{Span\left\{\frac{d}{dx}\vartheta_{u,v}(0):u\in \mathbb{R}_{+}^{m},v\in \mathbb{R}_{+}^{n}\right\}\right\}=mn,
\end{align*}
hence $\dim \mathcal{N}_{\mathcal{V}_{rn}}(E_{k})=0$, which shows that $E_{k}$ is singular.

Combine the above conclusions, we can get $\mathcal{V}_{rn}$ is an irreducible real algebraic variety with dimension $(m+n)r-r^2$, and $\mathcal{V}_{rn}^{ns}$ is
a $\mathbb{C}^\infty$-manifold of dimension $(m+n)r-r^2$. $\quad \Box$

\vspace{3mm}
\noindent
{\bf Proof of Theorem \ref{theorem2}}: Suppose that the angle $\alpha(A)$ of $A\in \mathcal{M}_{r}\cap \mathcal{M}_{n}$ is well defined, then $A$ is tangential if  $\alpha(A)=0$,
and is nontangetial if $\alpha(A)>0$.   Moreover, $A\in \mathcal{M}_{r}\cap \mathcal{M}_{n}$ is nontangential if and only if
\begin{align}
\mathcal{T}_{\mathcal{M}_{r}}(A)\cap \mathcal{T}_{\mathcal{M}_{n}}(A)=\mathcal{T}_{\mathcal{M}_{r}\cap\mathcal{M}_{n}}(A),
 \end{align}
 where $\mathcal{T}_{\mathcal{M}_{r}}(A)$ and $\mathcal{T}_{\mathcal{M}_{n}}(A)$ denote the differential geometry tangent spaces of manifolds $\mathcal{M}_{r}$ and $\mathcal{M}_{n}$ at point $A,$ respectively. Denote  $\mathcal{M}_{rn}^{nt}\subseteq \mathcal{M}_{rn}$ as the set of all  nontangential points of $\mathcal{M}_{r}\cap \mathcal{M}_{n}$.  Recall that the tangent space of $\mathcal{M}_{r}$ at $A=U_{m\times r}\cdot \Sigma_{r\times r} \cdot V_{n\times r}^{T}$ can be expressed as
 \begin{align*}
 {\cal T}_{{\cal M}_{r}}(A)
=\left\{[U,U_{\bot}]\left(
                                \begin{array}{cc}
                                  \mathbb{R}^{r\times r} & \mathbb{R}^{r\times (n-r)} \\
                                  \mathbb{R}^{(m-r)\times r} & 0^{(m-r)\times (n-r)} \\
                                \end{array}
                              \right)[V,V_{\bot}]^{T}\right\},
\end{align*}
where $U_{\bot}\in \mathbb{R}^{m\times (m-r)}$ and $V_{\bot}\in \mathbb{R}^{n\times (n-r)}$
stand for the orthogonal completions of $U$ and $V,$ respectively.
%In particular, $\mathcal{M}_{n}$ is a manifold with boundary so we need to introduce the tangent space of the points on the boundary. Assume that $P$ is a boundary point of $\mathcal{M}_{n}$, there are many method to define the tangent space to $\mathcal{M}_{n}$ at $P$, the standard choice is to define $T_{\mathcal{M}_{n}}(p)$ to be an mn-dimensional vector space. This may or may not seem like the most geometrically intuitive choice, but it has the advantage of making most of the definitions of geometric objects on a manifold with boundary look exactly the same as those on a manifold. Some details can be found in \cite{absil2009optimization,lee2013smooth}. After that i
In addition, it is not difficult  to derive $\mathcal{T}_{\mathcal{M}_{n}}(A)=Span\{E_{ij}\}$ with $E_{ij}=1$ or $0$,
for all $i=1,...,m$ and $j=1,...,n$.  Thus
\begin{align}\label{eq2}
\dim(\mathcal{T}_{\mathcal{M}_{r}}(A)\cap \mathcal{T}_{\mathcal{M}_{n}}(A))\leq (m+n)r-r^2.
\end{align}
Recall Theorem \ref{th1}, $\mathcal{M}_{rn}$ is a smooth manifold with dimension $(m+n)r-r^2$, which is homeomorphism to its tangent space $\mathcal{T}_{\mathcal{M}_{rn}}(A)$. Thus,
 \begin{align}\label{eq2}
\dim(\mathcal{T}_{\mathcal{M}_{r}}(A)\cap \mathcal{T}_{\mathcal{M}_{n}}(A))\leq \dim(\mathcal{T}_{\mathcal{M}_{rn}}(A)),
\end{align}
which is sufficient to say $\mathcal{M}_{rn}^{nt}$ is not empty.$\quad \Box$

Note that  $\mathcal{V}_{rn}$ is an irreducible variety with dimension $(m+n)r-r^2$, and $\mathcal{M}_{rn}^{nt}$ is not empty, then it follows from Theorem 6.6 in \cite{andersson2013alternating} that
$\mathcal{V}_{rn}\setminus \mathcal{M}_{rn}^{nt}$ is a real algebraic variety of dimension strictly less than $(m+n)r-r^2$. This result tell us the majority of points are nontangential if one is. Moreover,  \eqref{eq2} is satisfied, then by Theorem 5.1 in \cite{andersson2013alternating}, we can get Theorem \ref{thm_convergence}.

\end{document}